\documentclass[11pt]{article}

\usepackage{amsfonts}
\usepackage{amssymb,amsmath}
\usepackage{latexsym}
\usepackage{graphics}
\usepackage{epic}
\usepackage{epsfig}
\usepackage{subfigure}

\def\1{{\bf 1}}

\def\1{{\bf 1}}

\def\Section{\setcounter{equation}{0}\section}
\normalbaselines
\begin{document}
%\title{A Particle System describing Migration, Mutation and Selection
 % and some large population approximations}
% \title{Stochastic Individual Processes and Deterministic
%   Approximations for Spatially Structured Populations with Adaptive
%   Evolution}
\title{Invasion and adaptive evolution for individual-based spatially
  structured populations}
% \title{Nonlocal and local
%   interaction in stochastic discrete spatially structured populations
%   with adaptive evolution. Large population asymptotics}
 \author{Nicolas Champagnat$^1$,
  Sylvie M\'el\'eard$^{2,3}$}

\footnotetext[1] {Weierstrass Institute for Applied Analysis and
  Stochastics, Mohrenstrasse 39, Berlin,
champagn@wias-berlin.de}
 \footnotetext[2]{ Universit\'e Paris 10,
MODAL'X, 200 av.\ de la R\'epublique, 92000 Nanterre,
sylvie.meleard@u-paris10.fr}
 \footnotetext[3]{Fonctionnement et \'evolution des syst\`emes \'ecologiques, UMR 7625,
  ENS, 46 rue d'Ulm, 75005 Paris}

\maketitle

\setlength{\textwidth}{15cm}
\setlength{\textheight}{22cm}
\setlength{\oddsidemargin}{.5cm}
\setlength{\evensidemargin}{-.5cm}
\setlength{\topmargin}{-.5cm}
\setlength{\abovedisplayskip}{3mm}
\setlength{\belowdisplayskip}{3mm}
\setlength{\abovedisplayshortskip}{3mm}
\setlength{\belowdisplayshortskip}{3mm}
\raggedbottom
\parskip=1.5mm
\def\Section{\setcounter{equation}{0}\section}
\def\theequation {\arabic{section}.\arabic{equation}}
\def\be{\begin{eqnarray}}
\def\ee{\end{eqnarray}}
\def\ben{\begin{eqnarray*}}
\def\een{\end{eqnarray*}}
\def\numero{\refstepcounter{equation} (\theequation)}

\def\pit{\mathbb{P}}
\def\qit{\mathbb{Q}}
\def\rit{\mathbb{R}}
\def\r{\rit^2}
\def\nit{\mathbb{N}}
\def\dit{\mathbb{D}}
\def\eit{\mathbb{E}}
\def\One{\hbox{\it 1\hskip -3pt I}}
\def\e{\varepsilon}
\def\ala{\nonumber\\}

\newtheorem{thm}{Theorem}[section]
\newtheorem{lem}[thm]{Lemma}
\newtheorem{cor}[thm]{Corollary}
\newtheorem{defi}[thm]{Definition}
\newenvironment{definition}{\begin{defi} \rm}{\end{defi}}
\newtheorem{prop}[thm]{Proposition}
\newtheorem{rema}[thm]{Remark}
\newenvironment{remark}{\begin{rema} \rm}{\end{rema}}
\newenvironment{proof}{\noindent {\bf Proof}.}{\hfill$\Box$\\[-2mm]}
\newenvironment{display} {$$\displaylines{ } { \numero}$$}

\begin{abstract}

The interplay between space and evolution is an important issue in
population dynamics, that is in particular crucial in the
emergence of polymorphism and spatial patterns. Recently,
biological studies suggest that invasion and evolution are closely
related. Here we model the interplay between space and evolution
starting with an individual-based approach and show the important
role  of parameter scalings on clustering and invasion. We
consider a stochastic discrete model with birth, death,
competition, mutation and spatial diffusion, where all the
parameters may depend both on the position and on the trait of
individuals. The spatial motion is driven by a reflected diffusion
in a bounded domain. The interaction is modelled as a trait
competition between individuals within a given spatial interaction
range. First, we give an algorithmic construction of the process.
Next, we obtain large population approximations, as weak solutions
of nonlinear reaction-diffusion equations with Neumann's boundary
conditions. As the spatial interaction range is fixed, the
nonlinearity is nonlocal.  Then, we make the interaction range
decrease to zero and prove the convergence to spatially localized
nonlinear  reaction-diffusion equations, with Neumann's boundary
conditions. Finally, simulations based on the microscopic
individual-based model are given, illustrating the strong effects
of the spatial interaction range on the emergence of spatial and
phenotypic diversity (clustering and polymorphism) and on the
interplay between invasion and evolution. The simulations focus on
the qualitative differences between local and nonlocal
interactions.
\end{abstract}
%\vskip 1,5cm

 \vskip 5cm \pagebreak
\def\Section{\setcounter{equation}{0}\section}
\setlength{\normalbaselineskip}{17pt} \normalbaselines

\setcounter{equation}{0}
\section{Introduction}
\label{sec:intro}

The spatial aspect of population dynamics is an important
ecological issue that has been extensively studied (Murray
\cite{Murray89}, Durrett and Levin \cite{DL94b}, Tilman and
Kareiva \cite{TK96}, McGlade \cite{McG99}, Dieckmann et al.
\cite{DLM00}). It is in particular crucial in environmental
problems, such as spatial invasions and epidemics (Mollison
\cite{Mo77}, Murray~\cite{Murray89}, Rand et al. \cite{Rand-al95},
Tilman and Kareiva \cite{TK96}, Lewis and Pacala \cite{LP00}), and
clustering or agglomeration of the population, i.e.\ the
organization as isolated patches (Hassel and May \cite{HM74},
Hassel and Pacala \cite{HP90}, Niwa \cite{Ni94}, Flierl et
al.~\cite{FGLO99}, Young et al.~\cite{YRS01}).  The combination of
space and phenotype is also known for a long time to have
important effects on population dynamics (Mayr~\cite{Mayr63},
Endler~\cite{Endler77}).  In particular, it can strongly favor the
coexistence of several types of individuals and the emergence and
stability of polymorphism (Durrett and Levin~\cite{DL94},
Dieckmann and Doebeli~\cite{DD99}).  More recently, several
biological studies (Thomas et al. \cite{Tal01}, Phillips et al.
\cite{Pal06}) observed that  classical models could  underestimate
the invasion speed and suggested that evolution and invasion are
closely related. Namely, the evolution of morphology can have
strong impact on the expansion of invading species, such as insect
species (\cite{Tal01}) or cane toads (\cite{Pal06}). In this
context, the study of space-related traits, such as dispersal
speed (Prévost \cite{Prevost:04}, Desvillettes et al.
\cite{Desville:04}), or sensibility to heterogeneously distributed
resources (Bolker and Pacala~\cite{BP:99}, Grant and
Grant~\cite{GG02}), is fundamental.
% Conversely, the physical space can have strong effects on the
% phenotypic distribution (becs d'oiseaux).

In this paper, we propose and construct stochastic and
deterministic population models describing the interplay between
evolution and spatial structure. We show how helpful these models
can reveal to understand and predict several specific behaviors
concerning  clustering and invasion.

We study the dynamics of a spatially structured asexual population
with adaptive evolution, in which individuals can move, reproduce
with possible phenotypic mutation, or die of natural death or from
the competition between individuals. The individuals are
characterized both by their position and by one or several
phenotypical adaptive traits (such as body size, rate of food
intake, age at maturity or dispersal speed). The interaction is
modelled as a trait competition between individuals in some
spatial range. Our approach is based on a stochastic microscopic
description of these individuals' mechanisms, involving both space
and traits. This approach has already been developed in simpler
ecological contexts. For populations undergoing dispersal, Bolker
and Pacala \cite{Bolker:97, BP:99} and Dieckmann and Law
\cite{DL00}, offered the first microscopic heuristics and
simulations. Their individual-based model has been rigourously
constructed in Fournier and M\'el\'eard \cite{Fournier:04}. If one
thinks of the dispersion in the physical space as a mutation in a
trait space, this model translates into an evolutionary model. The
generalization to adaptive population with general mutation and
competition phenomena is achieved by Champagnat, Ferri\`ere and
M\'el\'eard \cite{Champagnat:05, Champagnat:06}.  In these papers,
different large population deterministic or stochastic
approximations have been obtained, depending on several scalings
on the microscopic parameters.

%   We construct a
% stochastic individual-based point process describing this
% population, from which we deduce large population approximations,
% with nonlocal nonlinearity when the spatial interaction range is
% fixed. Next, we make the interaction range decrease to zero and
% obtain some local nonlinear approximation, solving as corollary
% some nonlinear integro-differential equations. Finally, we give
% some simulations  illustrating the effect of the spatial
% interaction range in the emergence of phenotypic diversity
% (polymorphism) and  of spatial diversity (clustering).
%    The simulations also show the difference of qualitative behavior
%   between local and nonlocal interactions.

The basic mechanisms of the population dynamics we consider combine
spatial motion and evolutionary dynamics (Section~\ref{sec:model}).
The birth, mutation and death parameters of each individual depend on
its position and trait. An offspring, appearing at the same position
as its progenitor, usually inherits the trait value of the latter,
except when a mutation causes the offspring to take an instantaneous
mutation step at birth to new trait values.  As soon as it is alive,
an individual moves in the spatial domain according to a reflected
diffusion process. Moreover, each individual dies because of natural
death or is eliminated in the competition (selecting the fittest
traits) between individuals living in a given spatial range
$\delta>0$.

Section~\ref{sec:mart} starts with the algorithmic construction of a
stochastic Markov point process whose generator captures the
individual migration and ecological mechanisms in the population. Then
the existence of this measure-valued process and its martingale
properties are proved under some moment condition on the initial data.

Next (Section~\ref{sec:fixed-delta}), we study approximations of
this model based on large-population limits. We consider a large
number $N$ of individuals at initial time and assume that a fixed
amount of available resources has to be partitioned between
individuals.  When $N$ tends to infinity, the conveniently
normalized point process converges to a deterministic finite
measure, solution of a nonlinear nonlocal integro-differential
equation with Neumann's boundary conditions, parameterized by the
spatial range.  The proof is based on the martingale properties of
the process and on limit theorems for measure-valued jump
processes. We moreover prove that for sufficiently smooth and
non-degenerate diffusion coefficients, assuming that the initial
condition has a density, the limiting measure has at each time a
density with respect to the Lebesgue measure. That is due to the
regularizing effect of the reflected diffusion process.  The proof
mainly uses analytic tools, and is based on the mild formulation
of the limiting nonlinear equation.

In Section~\ref{sec:delta=0}, we study the behavior of this density
function as the interaction range tends to $0$. We show its
convergence to the solution of a spatially local nonlinear
integro-differential equation with Neumann's boundary conditions. This
equation has been introduced and studied in Pr\'evost
\cite{Prevost:04} in an analytic point of view (see also Desvillettes
et al. \cite{Desville:04}).
% Our microscopic construction justifies the biological interest of
% these local or nonlocal nonlinear models and implies existence results
% under minimal assumptions.
In this spatially local case, numerical simulations by finite
element methods are given and show the influence of diffusion and
mutation parameters on the invasion of the domain by the
population.

In Section~\ref{sec:simul}, we give simulations of the microscopic
process illustrating the time-dependent interplay between space and
adaptation. We address the effect of the population size, and the
crucial role of the interaction range with respect to spatial
organization (clustering) and polymorphism. We focus on the
qualitative differences between nonlocal and local interactions.  In a
first example, we show that, when migrations and mutations are not too
strong, a large interaction range induces a spatial organization of
the population as a finite set of isolated clusters, as assumed in
classical metapopulation models (\cite{DL94}).  Such a spatial
organization is related to the ecological notion of ``niches''
(different types of individuals settle different regions of space,
Roughgarden~\cite{Rough72}).  Conversely, for sufficiently small
interaction range, the clustering phenomenon is no more observed.
Next, we propose another example where a similar phase transition
occurs for spatial clustering and in which the critical interaction
range can be identified. In our last example, we investigate a model
describing the invasion of a species with evolving dispersal speed (as
in \cite{Desville:04}). The diffusion coefficient and the trait are
assumed to be proportional and a triangular invasion pattern is
observed, indicating that the invasion front is
composed of faster individuals (\cite{Pal06}).\\

\noindent{\bf Notation}\\
The individuals live in the closure of a bounded domain ${\cal X}$ of
$\rit^d$ of class $C^{3}$ and their trait values belong
to a compact set ${\cal U}$ of $\rit^q$.\\
- For $x\in\partial {\cal X}$, we denote by $\:n(x)\:$ the outward
normal to the boundary $\partial {\cal X}$ at point $\:x$.\\
- For a sufficiently smooth function $\: f\:$ and $(x,u)\in\partial
{\cal X}\times {\cal U}$, we denote by $\:\partial_nf(x,u)\:$ the
scalar product $\:\nabla_x
f(x,u)\cdot n(x)$.\\
- We denote by $\:C^{2,b}_0\:$ the space of measurable functions
$f(x,u)$ of class $C^{2}$ in $x$ and bounded in $u$ satisfying
$\partial_nf(x,u)=0$ for all $(x,u)\in
\partial {\cal X}\times {\cal U}$
and by $\:C^{2,0}_0\:$ the subspace of functions $f(x,u)$ which are
moreover continuous in $u$.\\
- For each $p\geq 1$, the $L^p$-norm on $\bar{\cal X}\times {\cal U}$
 is denoted by $\|\cdot\|_{p}$.\\
% - If $E$ is a Polish space, the space of probability measures on
% $E$ is denoted by ${\cal P}(E)$.\\
- We denote by $M_F(\bar{\cal X}\times {\cal U})$ the set of finite
measures on $\bar{\cal X}\times {\cal U}$, endowed by the weak
topology, and by ${\cal M}$ the subset of $M_F(\bar{\cal X}\times
{\cal U})$ composed of all finite point measures, that is
\begin{equation*}
  {\cal M} = \left\{ \sum_{i=1}^n \delta_{(x^i,u^i)} , \; n \in \nit,\  x^1,\ldots,x^n
    \in \bar{\cal X} \ ,\ u^1,...,u^n \in {\cal U}\right\}
\end{equation*}
where $\delta_{(x,u)}$ denotes the Dirac measure at $(x,u)$. (If
$n=0$, one obtains by extension the null measure).  For any $\nu \in
M_F(\bar{\cal X}\times {\cal U})$ and for any measurable function $f$
on $\bar{\cal X}\times {\cal U}$, we write indifferently $\left< \nu,
  f \right>$ or $\int_{\bar{\cal X}\times {\cal U}} f
d\nu$. If $\nu=\sum_{i=1}^n \delta_{(x^i,u^i)}$, then
$\langle\nu,f\rangle=\sum_{i=1}^nf(x^i,u^i)$ .\\
- We denote by $\dit ([0,\infty), M_F(\bar{\cal X}\times {\cal U}))$
the Skorohod space of left limited and right continuous functions from
$\rit_+$ to
$M_F(\bar{\cal X}\times {\cal U})$, endowed with the Skorohod topology.\\
- The constant $C$ will be a constant which can change from line to
line.

\begin{rema}
  \label{densite} Let us remark that the space of $C^2(\bar{\cal
    X})$-functions with a vanishing normal derivative is dense, for
  the uniform norm, in $C(\bar{\cal X})$. Indeed, let us consider the
  Cauchy problem for the parabolic differential equation ${\partial
    u\over \partial t}(t,x)= \Delta u(t,x)\ ;\ t>0\ ; \ x\in {\cal X}$
  with the boundary condition ${\partial u\over \partial n}(t,x)= 0\
  ;\ t>0\ ; \ x\in \partial{\cal X}$. Since ${\cal X}$ is of class
  $C^{3}$, we may apply Sato-Ueno \cite{Sato:65} Theorem 2.1. There
  exists a smooth fundamental solution $q(t,x,y)$ to this problem and
  each $f\in C(\bar{\cal X})$ is the uniform limit of the sequence
  $\:\int_{\bar{\cal X}}q(t,x,y)f(y)dy\:$ of $C^2(\bar{\cal
    X})$-functions with vanishing normal derivative, as $t$ tends to
  $0$.

  We easily extend this result and show that the space $C^{2,0}_0$ is
  dense in the space of continuous functions on $\bar{\cal X}\times
  {\cal U}$.
\end{rema}

\setcounter{equation}{0}
\section{The model}
\label{sec:model}

Let us now describe the evolutionary process we are interested in. The
population will be described at any time by a finite point measure
$\nu\in {\cal M}$. Each individual, characterized by its position and
trait $(x,u)$, may move, give birth or die, as described below.
\begin{enumerate}
\item The {\bf migration} is described as a diffusion process normally
  reflected at the boundary of the domain ${\cal X}$. Biologists
  usually assume that the random behavior is isotropic, so the
  diffusion matrix is chosen with the form $m(x,u) \mbox{Id}$ ($\mbox{Id}$ is the
  identity matrix on $\rit^d$) and the nonnegative coefficient
  $m(x,u)$ (depending on the position $x$ and the trait value $u$), is the
  diffusion coefficient. We moreover model the environment heterogeneity
  (resources, topography, external effects,\ldots) by a drift term
  driven by a $\rit^d$-vector $b(x,u)$.
\item {\bf Births and mutations}.  We consider a population with
  asexual reproduction.  An individual with position $x$ and trait
  $u$ can give birth either to a clonal child at rate $\lambda(x,u)$,
  or to a mutant with trait $v$ at rate $M(x,u,v)$, both at position
  $x$.
\item The {\bf death} rate $\:\mu\:$ of an individual depends on its
  position $x$ and trait $u$ and on the spatial and phenotypic
  competition with the individuals located around $x$.
  Let us call $\delta>0$ the \textbf{range} of this spatial interaction.

  For a population $\nu=\sum_{i=1}^n \delta_{(x^i,u^i)}\in {\cal M}$,
  the death rate is given by
  \begin{align*}
    \mu(x,u,I^\delta W\star \nu(x,u))&=\mu\biggl(x,u,\int_{{\cal
        X}\times {\cal U}}I^\delta(x-y)W(u-v)\nu(dy,dv)\biggr)\\
    &=\mu\biggl(x,u,\sum_{i=1}^nI^\delta(x-x^i) W(u-u^i)\biggr).
  \end{align*}
  The function $\mu(x,u,r)$ is assumed to be measurable on ${\cal
    X}\times {\cal U}\times \rit$.
\end{enumerate}

This interaction assumes that spatial and phenotypic interactions
are independent, which is realistic in many biological situations.
One could of course consider a more complicated interaction. Since
our ultimate goal is to make the spatial
interaction range go to zero, we have chosen this particular form.\\\\
{\bf Hypotheses (H):}\\ {\it 1) The coefficients $m(x,u)$ and $b(x,u)$
  depend Lipschitz continuously on the position and measurably on the
  trait, and there exist constants $m^\star>0$ and $b^*>0$ such that
  for all $(x,u)\in \bar{\cal X}\times {\cal U}$
  \begin{equation}
    \label{migration}
    \begin{gathered} 0\leq m(x,u)\leq m^\star \\
      |b(x,u)|\leq b^*. \end{gathered}
  \end{equation}
  2) It is natural from a biological point of view to assume that all
  birth rates are bounded. There exists $\lambda^*$ such that
  \begin{equation}
    \label{birth} 0\leq \lambda(x,u)\leq \lambda^*,\quad \forall
    (x,u)\in \bar{\cal X}\times {\cal U}.
  \end{equation}
  The kernel $M$ is nonnegative and symmetric in $(u,v)$ for each
  $x\in \bar{\cal X}$ and
  \begin{equation}
    \label{mutation} \sup_{x\in \bar{\cal X},u\in
      {\cal U}}M(x,u,v)=M^*(v) \in L^1({\cal U}).
  \end{equation}
  3) There exists a positive constant $\mu^*$ such that
  \begin{equation}
    \label{death} \forall (x,u,r) \in {\cal X}\times {\cal U} \times \rit,
    \quad 0\leq \mu(x,u,r)\leq \mu^*(1+|r|)
  \end{equation}
  4) For each $\delta>0$, the spatial kernel $I^\delta$ is nonnegative
  and bounded and
  for each $x\in \bar{\cal X}$,
 $$\int_{{\cal
      X}}I^\delta(x-y)dy=1.$$
5) The competition kernel $W$ is
  nonnegative and bounded on $\rit^q$.}\\

Let us remark that if $I^\delta$ is proportional to $1_{\{|x|\leq
  \delta\}}$, then (H-4) means that the interaction is proportional to
the surface in ${\cal X}$ around $x$. This is a natural biological
assumption, especially if $x$ lies on the boundary of ${\cal X}$. We
will later assume that the measure $I^{\delta}(y)dy$ weakly converges
to the Dirac measure $\delta_0$ as $\delta$ tends to $0$.
\\

Hypotheses (H) will be assumed in all the sequel. They imply in
particular that for each $\nu\in M_F(\bar{\cal X}\times {\cal U})$
and each $(x,u )\in\bar{\cal X}\times {\cal U}$,
\begin{equation}
  \label{cdelta}
  \mu(x,u,I^\delta W\star\nu(x,u))\leq \mu^*(1+\|I^\delta
  W\|_\infty\langle\nu,1\rangle)
\end{equation}
which yields
\begin{align}
  &\mu(x,u,I^\delta
  W\star\nu(x,u))  +\lambda(x,u)+\int_{{\cal U}}M(x,u,v)dv\notag \\
  &\hskip 2cm\leq \mu^*(1+\|I^\delta
  W\|_\infty\langle\nu,1\rangle)+\lambda^*+ \|M^*\|_{1}
  \leq C_\delta(\langle\nu,1\rangle+1) \label{cdelta1}
\end{align}
and the total jump rate for a population $\nu$ is bounded by \be
\label{totjump}
C_\delta\langle\nu,1\rangle(\langle\nu,1\rangle+1).\ee

 We are interested in the evolution of the stochastic point process $(\nu_t)$,
taking its values in ${\cal M}$ and describing the evolution of
the population at each time $t$. We define
\begin{equation*}
  \nu_t = \sum_{i=1}^{N_t} \delta_{(X^i_t, U^i_t)}\ ,
\end{equation*}
$N_t \in {\nit}$  standing for the number of living individuals at
time $t$,  $X^1_t,...,X^{N_t}_t$  describing their positions (in
$\bar{\cal
  X}$)  and $U^1_t,...,U^{N_t}_t$  their trait values (in ${\cal U}$).

The dynamics of the population can be roughly summarized as follows.
The initial population is characterized by a measure $\nu_0 \in {\cal
  M}$ at time $t=0$, and any individual located at $x \in \bar{\cal
  X}$ with trait $u$ at time $t$ has four independent exponential
clocks: a ``clonal reproduction'' clock with parameter $\lambda(x,u)$,
a ``mutant reproduction'' clock with parameter $M(x,u,v)$, and a
``mortality'' clock with parameter $\mu(x,u,
\sum_{j=1}^{N_t}I^\delta(x-X^j_t)W(u-U^j_t))$. If the ``mortality''
clock of an individual rings, then this individual disappears; if the
``clonal reproduction'' clock of an individual rings, then it produces
at the same location an individual with the same trait as itself; if
the ``mutant reproduction'' clock of an individual rings, then it
produces at the same location an individual with characteristics
$(x,v)$.

The living individuals evolve in the domain, according to
diffusion processes with diffusion coefficient $m(x,u)$ and drift
$b(x,u)$, normally reflected at the boundary of ${\cal X}$.

The measure-valued process $(\nu_t)_{t\geq 0}$ is a Markov process
whose infinitesimal generator $L$ captures this dynamics.  This
generator is the sum of a jump part $L_1$ corresponding to the
phenotypic evolution and of a diffusion part $L_2$. The generator
$L_1$ is defined for bounded and measurable functions $\phi$ from
${\cal M}$ into $\rit$ and for $\ \nu=\sum_{i=1}^n
\delta_{(x^i,u^i)}\ $ by
\begin{align}
  L_1\phi{}(\nu)&=\sum_{i=1}^{\langle\nu,1\rangle}\lambda(x^i,u^i)(\phi{}(\nu+\delta{}_{(x^i,u^i)})-\phi{}(\nu))\notag\\
  &+\int_{\cal
    U}\sum_{i=1}^{\langle\nu,1\rangle}(\phi{}(\nu+\delta{}_{(x^i,v)})-\phi{}(\nu))
  M(x^i,u^i,v)dv\notag \\
  &+\sum_{i=1}^{\langle\nu,1\rangle}(\phi(\nu-\delta_{(x^i,u^i)})-\phi{}(\nu))
  \mu{}(x^i,u^i,I^\delta W\star\nu(x^i,u^i)). \label{gensaut1}
\end{align}

A standard class of cylindrical functions generating the set of
bounded and measurable functions from $M_F(\bar{\cal X}\times {\cal
  U})$ into $\rit$ is the class of functions
\begin{equation}
  \label{cylind}
  F_f(\nu)=F(\langle \nu,f\rangle),
\end{equation}
for bounded and measurable functions $F$ and $f$.

For such functions $F_f$, with $F\in C^2_b(\rit)$ and $f\in
C^{2,0}_0$, the diffusive part $L_2$ of the generator can easily be
deduced from It\^o's formula. Its form is similar to the one obtained in
the whole space for branching diffusing processes (cf.  Roelly-Rouault
\cite{Roelly:90}) and is given by
\begin{equation}
  \label{gendiff} L_2 F_f(\nu)=\langle\nu, m\Delta_x f+b.\nabla_x
  f\rangle F'(\langle\nu,f\rangle)+\langle\nu, m|\nabla_x f|^2\rangle
  F''(\langle\nu,f\rangle).
\end{equation}
Hence,
\begin{align}
  LF_f(\nu)&=L_1F_f(\nu)+L_2F_f(\nu)\notag\\&=\int_{\bar{\cal X}\times
    {\cal U}} \bigg\{\lambda(x,u)\big(F(\langle\nu,f\rangle+f(x,u))-
  F(\langle\nu,f\rangle)\big)\notag\\ &\qquad\qquad\ +\int_U
  \big(F(\langle\nu,f\rangle+f(x,v))-
  F(\langle\nu,f\rangle)\big)M(x,u,v)dv\notag\\ &\qquad\qquad\
  +\mu(x,u,I^\delta W\star \nu(x,u))\big(F(\langle\nu,f\rangle-f(x,u))-
  F(\langle\nu,f\rangle)\big)\notag\\ &\qquad\qquad\ +\big(
  m(x,u)\Delta_xf(x,u)+b(x,u).\nabla_x f(x,u)\big) F'(\langle\nu,f\rangle))\notag\\
  &\qquad\qquad\ + m(x,u)|\nabla_xf(x,u)|^2 F''(\langle\nu,f\rangle)\bigg\}\
  \nu(dx,du)\label{generator}
\end{align}

\setcounter{equation}{0}
\section{Construction of the particle system and martingale
  properties}
\label{sec:mart}

In this section, we construct a Markov process on the path space $\dit
([0,\infty), M_F(\bar{\cal X}\times {\cal U}))$ with infinitesimal
generator $L$. Then we prove some martingale properties satisfied by
this process, which are the key point to obtain large population
approximations.

Let us firstly present an iterative construction of the process,
which gives an effective simulation algorithm, if combined with a
diffusion simulation step such as an Euler scheme for reflected
diffusions (see L\'epingle \cite{Le95}, Gobet \cite{Gobet:01} and
Section~\ref{sec:simul}).

The initial number of individuals is equal to some natural integer
$N\in\nit^*$ and the vector of random variables
$(X_0,U_0)=(X^{i}_0, U^i_0)_{1\leq i\leq N}\in (\bar{\cal X}\times
{\cal U})^N$ denotes the position and trait values of these
individuals. More generally, we denote by $N_t$ the number of
individuals at time $t$ and by $(X_t,U_t)$ the vector of their
positions and traits. Let us introduce the following sequences of
independent random variables, independent of $(X_0,U_0)$.

-  $(B^{j,k})_{k,j\in \nit^*}$ are  $d$-dimensional Brownian
motions,

- $(\theta_k)_k$ are uniform random variables on $[0,1]$,

- $(V_k)_k$ take values in ${\cal U}$ with law
 ${M^*(v)\over
  \|M^*\|_{1}}dv$,

- $(\tau_k)_k$ are exponential random variables with law $C_\delta
e^{-C_\delta t}1_{t\geq 0}$. (The constant $C_{\delta}$ is defined in
(\ref{cdelta1})).

The system is obtained inductively for $k\geq 1$ as described below.
We set $T_0=0$ and $N_0=N$. Assume that
$(T_{k-1},N_{k-1},X_{T_{k-1}},U_{T_{k-1}})$ are given. If $N_{k-1}=0$,
then $\nu_t=0$ for all $t\geq T_{k-1}$. If not, let

\begin{itemize}
\item $T_{k}=T_{k-1}+\frac{\tau_k}{N_{k-1}(N_{k-1}+1)}$. Notice that
  ${\tau_k\over N_{k-1}(N_{k-1}+1)}$ represents the time between
  possible jumps for $N_{k-1}$ individuals and that $C_\delta
  N_{k-1}(N_{k-1}+1)$ gives an upper-bound on the total jump rate
  for a population with $N_{k-1}$ individuals, as seen in~(\ref{totjump}).

\item On the time-interval $[T_{k-1},T_k)$, the number of
particles
  remains equal to $N_{k-1}$, their trait values are equal to
  $U^j_{T_{k-1}}, 1\leq j\leq N_{k-1}$ and  their positions
  $(X^{j}_t,\;1\leq j\leq N_{k-1})$ evolve according to the following
  stochastic differential equation with normal reflection : $\forall
  t\in [T_{k-1},T_k]$,
  \begin{align}
    &X^{j}_t\in \bar{\cal X},\notag\\
    &X^{j}_t=X^{j}_{T_{k-1}}+\int_{T_{k-1}}^t
    \sqrt{2m(X^j_{s},U^j_{T_{k-1}})}dB^{j,k}_s
    +\int_{T_{k-1}}^tb(X^j_s,U^j_{T_{k-1}})ds -k^{j}_t\ , \notag\\
    &|k^{j}|_t=\int_{T_{k-1}}^t \1_{\{X^{j}_s\in
      \partial {\cal X}\}}d|k^{j}|_s\ ;\ k^{j}_t=\int_{T_{k-1}}^t
    n(X^{j}_s)d|k^{j}|_s.\label{sys.kps}
  \end{align}
\item At time $T_k$, one chooses at random an individual $I_k=i$
  uniformly among the $N_{k-1}$ individuals living during the
  time-interval $[T_{k-1},T_k)$. Its position and trait are
  $(X^i_{T_k},U^i_{T_{k-1}})$.
  % Remark that $C_{\delta}(N_{k-1}+1)$
  %gives an upper-bound on the total event rate for this individual.
  \begin{itemize}
  \item If $\:0\leq \theta_k\leq
    \frac{\mu(X^{i}_{T_{k}},U^{i}_{T_{k-1}},
      \sum_{j=1}^{N_{k-1}}I^\delta(X^i_{T_k}-X^j_{T_k})W(U^{i}_{T_{k-1}}-U^{j}_{T_{k-1}}))}
    {C_\delta(N_{k-1}+1)}=\theta_1^i(X_{T_{k}},U_{T_{k-1}})$, then the
    individual $\: i\: $ dies and $N_k=N_{k-1}-1$.
  \item If $\:\theta_1^i(X_{T_{k}},U_{T_{k-1}})< \theta_k\leq
    \theta_1^i(X_{T_{k}},U_{T_{k-1}})+{\lambda(X^{i}_{T_{k}},U^{i}_{T_{k-1}})
      \over C_\delta(N_{k-1}+1)}= \theta_2^i(X_{T_{k}},U_{T_{k-1}})$, then
    the individual $\: i\: $ gives birth to an offspring with
    characteristics $(X^i_{T_{k}},U^i_{T_{k-1}})$ and $N_k=N_{k-1}+1$.
  \item If $\:\theta_2^i(X_{T_{k}},U_{T_{k-1}})< \theta_k \leq
    \theta_2^i(X_{T_{k}},U_{T_{k-1}})+{M(X^{i}_{T_{k}},U^{i}_{T_{k-1}},V_k)\|M^*\|_1
      \over M^*(V_k)C_\delta(N_{k-1}+1)}=
    \theta_3^i(X_{T_{k}},U_{T_{k-1}},V_k)$, then the individual $\: i\: $
    gives birth to a mutant offspring with trait $V_k$ at the position
    $X^{i}_{T_{k}}$, and $N_k=N_{k-1}+1$.
  \item If $\: \theta_k> \theta_3^i(X_{T_{k}},U_{T_{k-1}},V_k)$, nothing happens
    and $N_k=N_{k-1}$.
  \end{itemize}
\end{itemize}

The total number $N_t$ of individuals at time $\,t\,$ is equal to
$N_t=\sum_{k\geq 0}1_{\{T_k\leq t<T_{k+1}\}}N_k$, and $ \nu_t=
\sum_{k\geq 0}1_{\{T_k\leq t<T_{k+1}\}}
\sum_{i=1}^{N_k}\delta_{(X^i_{t},U^i_{T_k})}=\sum_{i=1}^{N_t}\delta_{(X^i_t,U^i_t)}$.\\

This stochastic individual-based process $\nu$ can be rigorously
expressed as solution of a stochastic differential equation driven
by $d$-dimensional Brownian motions $(B^i)_{i\in \nit^*}$ and the
$\rit_+\times \nit\times[0,1]\times {\cal U}$-valued multivariate
point process
$$Q(dt,di,d\theta,dv)=\sum_{k\geq 1}\delta_{(T_{k},I_k,\theta_k,V_k)}(dt,di,d\theta,dv)$$
associated with the birth, mutation and death of individuals. We will
prove its existence on $\rit^+$, deduced from moment properties, and
develop some martingale properties that we will use below.

Let us consider $\nu_0\in {\cal M}$. For each $C^{2,0}_0$-function
$f$, we define the process $\langle \nu_t,f\rangle$ as solution of the
stochastic differential equation
\begin{align}
  &\langle
  \nu_t,f\rangle=\langle \nu_0,f\rangle+\int_0^t\langle
  \nu_r,m(x,u)\Delta_x f+b(x,u).\nabla_xf\rangle dr\notag\\ &+\int_0^t
  \sum_{i=1}^{\langle \nu_{r-},
    1\rangle}\sqrt{2m(X^i_{r},U^i_{r})}\nabla_xf(X^i_r,U^i_r)dB^i_r\notag\\
  &+\int_{[0,t]\times \nit \times [0,1]\times
    U^2}\bigg\{-f(X^i_r,U^i_r) \1_{\{\theta\leq \theta^i_1(X_r,U_r)\}} +
  f(X^i_r,U^i_r) \1_{\{\theta^i_1(X_r,U_r)<\theta\leq \theta^i_2(X_r,U_r)\}}\notag\\ &
  \phantom{+\int_{[0,t]\times \nit \times [0,1]\times
    U^2}\bigg\{\ }+ f(X^i_r,v) \1_{\{\theta^i_2(X_r,U_r)< \theta\leq
    \theta^i_3(X_r,U_r,v)\}}\bigg\}Q(dr,di,d\theta,dv), \label{sde}
\end{align}
where $\theta^i_1$, $\theta^i_2$ and $\theta^i_3$ have been defined previously.

By Remark \ref{densite}, the knowledge of $\langle \nu_t,f\rangle$ for
$f\in C^{2,0}_0$ is enough to characterized the finite measure-valued
process $\nu$.

We introduce the canonical filtration
$${\cal F}_t=\sigma\{\nu_0;\  B^j_r, j\in
\nit^*;\ Q([0,r]\times A ),\;A\in{\cal P}(\nit)\otimes{\cal
  B}([0,1]\times {\cal U}),\;r\leq t\},$$ where ${\cal B}([0,1]\times
{\cal U})$ is the Borel $\sigma$-field on $[0,1]\times{\cal U}$.

\begin{lem}
  \label{compens} The measure
  \begin{align*}
    q(dt,di,d\theta,dv)&= C_\delta
    \sum_{k\geq 0}1_{\{T_{k}< t\leq
      T_{k+1}\}}(N_k+1)\sum_{j=1}^{N_k}\delta_j(di) dtd\theta{M^*(v)\over
      \|M^*\|_1}dv \\ &=C_\delta
    (N_{t}+1)\sum_{j=1}^{N_{t}}\delta_j(di)dtd\theta{M^*(v)\over
      \|M^*\|_1}dv
  \end{align*}
  is the (predictable) compensator of the multivariate point process
  $Q$.
\end{lem}

\begin{proof}
  For $k\geq 0$, a regular version of the conditional law of
  $(T_{k+1},I_{k+1},\theta_{k+1},V_{k+1})$ with respect to $\sigma\{\nu_0,\;
  (B^j_.),\; j\in \nit^*,\;(T_p,I_p,\theta_p,V_p),\;1\leq p\leq k\}$ is
  given by the measure
  \begin{equation*}
    C_\delta (N_k+1)1_{\{T_k<t\}}e^{-C_\delta
      N_k(N_k+1)
      (t-T_k)}\;\sum_{j=1}^{N_k}\delta_j(di)dtd\theta{M^*(v)\over
      \|M^*\|_1}dv.
  \end{equation*}
  The conclusion is thus a consequence of \cite{Jacod:87} Theorem 1.33
  p.136.
\end{proof}

Using Lemma \ref{compens} and It\^o's formula, one can immediately
show that any solution $\nu$ of (\ref{sde}), such that
$E(\sup_{t\leq
  T}\langle\nu_t,1\rangle^2)<+\infty$, is a Markov process with
infinitesimal generator $L$ defined by (\ref{generator}).  Moreover,
we also deduce the following existence, moment and martingale
properties.

\begin{prop}
  \label{moment} 1) Assume Hypotheses (H) and that $E(\langle
  \nu_0,1\rangle)<+\infty$.

  Then $\ E(\sup_{t\leq T}\langle\nu_t,1\rangle)<+\infty\ $ for each
  $T>0$ and the process $\nu$ defined by (\ref{sde}) is well defined
  on $\rit^+$.

  \noindent 2) If furthermore for some $p\geq 1$,
  $E(\langle\nu_0,1\rangle^p)<+\infty$, then for each
  $T>0$
  \begin{equation*}
% \label{moment}
  E(\sup_{t\leq
    T}\langle\nu_t,1\rangle^p)<+\infty.
\end{equation*}
\end{prop}

\begin{proof} We firstly prove 2).  For each integer $k$, define
  $S_k=\inf\{t\geq 0,\: \langle\nu_t,1\rangle\geq k\}$. A simple
  computation using (\ref{sde}), and dropping the non-positive death
  terms, gives
  \begin{align*}
    E(\sup_{s\in [0,t\wedge
      S_k]}\langle\nu_s,1\rangle^p)&\leq
    E\left(\langle\nu_0,1\rangle^p+C \int_0^{t\wedge S_k}
      (1+\langle\nu_s,1\rangle^p)ds\right)\\ &\leq C\left(1+
      E\left(\int_0^t\langle\nu_{s\wedge S_k}
        ,1\rangle^pds\right)\right).
  \end{align*}
  Gronwall's lemma implies that for any $T>0$, there exists a constant
  $C$ independent of $k$, such that $E(\sup_{t\in[0,T\wedge
    S_k]}\langle\nu_t,1\rangle^p)\leq C$. One easily deduces that
  $S_k$ tends a.s.\ to infinity when $k$ tends to infinity and
  next, Fatou's lemma yields $E(\sup_{t\in
    [0,T]}\langle\nu_t,1\rangle^p)< +\infty$.

  Point 1) is a consequence of point 2). Indeed, one builds the
  solution $(\nu_t)_{t\geq 0}$ step by step. One only has to check
  that the sequence of jump instants $(T_k)_k$ goes to infinity a.s.\
  as $k$ tends to infinity. But this follows from $\: E(\sup_{t\leq
    T}\langle\nu_t,1\rangle)<+\infty\ $.
\end{proof}

The following martingale properties are the key point to study large
population approximations.

\begin{thm}
  \label{mpm} Assume Hypotheses (H) and that for some $p\geq 2$,
  $E(\langle\nu_0,1\rangle^p)<+\infty$.

  \noindent 1) Then, for $F$ and $f \in C^{2,0}_0$ such that for all
  $\nu\in{\cal M}$, $|F_f(\nu)|+|LF_f(\nu)|\leq
  C(1+\langle\nu,1\rangle^p)$, the process
  $$
  F_f(\nu_t)-F_f(\nu_0)-\int_0^t LF_f(\nu_s)ds
  $$
  is a c\`adl\`ag martingale starting from $0$.  It is in particular
  true for $F(y)=y^{p-1}$.

  \noindent 2) The process $Z^f$ defined for $f \in C^{2,0}_0$ by
  \begin{multline}
    Z^f_t=\langle \nu_t,f\rangle-\langle
    \nu_0,f\rangle- \int_0^t\int_{\bar{\cal X}\times
      U}\bigg\{m(x,u)\Delta_x f(x,u)+ b(x,u).\nabla_x f(x,u) \\
    +\big(\lambda(x,u) -\mu(x,u,I^\delta W\star\nu_s(x,u))\big)f(x,u)
    +\int_U f(x,v)M(x,u,v)dv \bigg\}\nu_s(dx,du) ds\label{mart1}
  \end{multline}
  is a c\`adl\`ag $L^2$-martingale starting from $0$ with predictable
  quadratic variation
  \begin{multline}
    \langle
    Z^f\rangle_t=\int_0^t\int_{\bar{\cal X}\times {\cal U}} \bigg\{2
    m(x,u)|\nabla_xf|^2+ \big(\lambda(x,u)+\mu(x,u,I^\delta
    W\star\nu_s(x,u))\big)f^2(x,u) \\ +\int_{\cal U} f^2(x,v)M(x,u,v)dv\bigg\}
    \nu_s(dx,du) ds  \label{quadra}
  \end{multline}
\end{thm}

\begin{proof}
  Point 1) is immediate. For point 2), we first assume that
  $E(\langle\nu_0,1\rangle^3)<+\infty$. Applying point 1) with
  $F(y)=y$ (or (\ref{sde}) and Lemma \ref{compens}) leads to $Z^f$.
  Then one applies 1) again with $F(y)=y^2$, and thus
  \begin{align}
    \langle \nu_t,f\rangle^2&-\langle \nu_0,f\rangle^2-
    \int_0^t\int_{\bar{\cal X}\times {\cal U}}\bigg\{2(m(x,u)\Delta_x f(x,u)+
    b(x,u).\nabla_x f(x,u))\langle\nu_s,f\rangle\notag\\
    &+2m(x,u)|\nabla_xf|^2+ \lambda(x,u)(2\langle\nu_s,f\rangle
    f(x,u)+f^2(x,u))\notag\\ & +\int_{\cal U}
    (2f(x,v)\langle\nu_s,f\rangle+f^2(x,v))M(x,u,v)dv \notag\\ &+\mu(x,u,I^\delta
    W\star\nu_s(x,u))(-2\langle\nu_s,f\rangle
    f(x,u)+f^2(x,u))\bigg\}\nu_s(dx,du)ds \label{carre}
  \end{align}
  is a c\`adl\`ag martingale. In another hand, It\^o's formula allows
  us to compute $\langle \nu_t,f\rangle^2$ from (\ref{mart1}): the
  process
  \begin{align}
    &\langle \nu_t,f\rangle^2 - \langle \nu_0,f\rangle^2 -
    \int_0^t\int_{\bar{\cal X}\times {\cal U}} \bigg\{2(m(x,u)\Delta_x f(x,u)+
    b(x,u).\nabla_x f(x,u))\langle\nu_s,f\rangle\notag\\ &+ 2(\lambda(x,u)-
    \mu(x,u,I^\delta W\star\nu_s(x,u)))\langle\nu_s,f\rangle
    f(x,u)\notag\\ & + \int_{\cal U}
    2f(x,v)\langle\nu_s,f\rangle M(x,u,v)dv\bigg\}\nu_s(dx,du) ds -
    \langle Z^f\rangle_t \label{ito}
  \end{align}
  is a c\`adl\`ag martingale. Comparing (\ref{carre}) and (\ref{ito})
  leads to (\ref{quadra}). The extension to the case where
  $E(\langle\nu_0,1\rangle^2)<+\infty$ is straightforward, noticing
  that $E(\langle Z^f\rangle_t)<+\infty$.
\end{proof}

\setcounter{equation}{0}
\section{Large population approximation for a fixed interaction range}
\label{sec:fixed-delta}

We are now interested in deterministic approximations of the
population point process when the size of the population increases. We
assume in this section that the interaction range $\delta>0$ is fixed.

Let us consider a sequence of initial measures
$(\nu^N_0)_{N\in\nit^*}$ belonging to ${\cal M}$. For each
$N\in\nit^*$, we keep all parameters $(m,b,\lambda,M)$ unchanged,
except the competition kernel. We assume that for each $N$,
\begin{equation}
  \label{sn} \mu_N(x,u,r)=\mu(x,u,{r\over N}).
\end{equation}
This assumption has a biological interpretation. In a case of fixed
amount of available global resources, a large system of individuals
may only exist if the biomass of each interacting individual scales as
${1\over N}$, which implies that the interaction effect between two
individuals scales as ${1\over N}$ as well. The parameter $N$ can also
be interpreted as scaling the resources available, so that the
renormalization of $\mu$ reflects the decrease of competition for
resources.

We assume that the sequence ${\nu^N_0\over N}$ converges, as $N$ tends
to infinity. The size $\langle\nu_0^N,1\rangle$ of the population is
then of order $N$ and will stay at this order (or at a smaller order)
during finite time-intervals, since birth rates are bounded. Hence,
our aim is to study the asymptotic behavior, as $N$ tends to infinity,
of the c\`adl\`ag process
\begin{equation}
  \label{mvp}
  \Lambda^N_t=\frac{1}{N}\sum_{i=1}^{N_t}\delta_{(X^{i}_t,U^i_t)}=\frac{1}{N}\nu^N_t,
\end{equation}
taking values in ${\cal M}^N=\{{1\over N}\nu, \nu\in {\cal M}\}$.

The process $(\Lambda^N_t)_{t\geq 0}$ is a Markov process with
generator $L_N=L_{N,1}+L_{N,2}$. An easy computation, for $F\in
C^2(\rit)$ and $f\in C^{2,0}_0$, gives that
\begin{equation}
  \label{gendiffap}
  L_{N,2}F_f(\nu)=\langle \nu,m(.)\Delta_xf+b(.).\nabla_x f\rangle
  F'(\langle\nu,f\rangle) +\langle \nu,{m(.)\over N}|\nabla_x
  f|^2\rangle F''(\langle\nu,f\rangle)
\end{equation}
and (using (\ref{sn}))
\begin{align}
  L_{N,1}F_f(\nu)&=N\int_{\bar{\cal X}\times {\cal
      U}}\bigg\{\lambda(x,u)\big(F(\langle \nu,f\rangle+{1\over
    N}f(x,u))-F(\langle \nu,f\rangle)\big)\notag\\ &+ \mu(x,u,I^\delta
  W\star\nu(x,u))\big(F(\langle \nu,f\rangle-{1\over N}f(x,u))-F(\langle
  \nu,f\rangle)\big)\notag\\ &+\int_{\cal U}\big(F(\langle \nu,f\rangle+{1\over
    N}f(x,v))-F(\langle \nu,f\rangle)\big)M(x,u,v)dv\bigg\}\nu(dx,du)\label{gensautap}
\end{align}

We deduce from Theorem \ref{mpm} the following martingale properties.

\begin{lem}
  \label{martn}
  Let $N\geq 1$ be fixed and assume that for some $p\geq 2$,
  $E \left( \left< \Lambda_0^N,1 \right>^p \right) <\infty$.
 For all $C^{2,0}_0$-function $f$, the process
  \begin{multline}
    Z^{N,f}_t = \left<\Lambda^N_t,f\right>-
    \left<\Lambda^N_0,f\right> -\int_0^t\int_{\bar{\cal X}\times {\cal U}}
    \bigg\{m(x,u)\Delta_x f(x,u) + b(x,u).\nabla_x f(x,u) \\ +
    \big(\lambda(x,u)
    -\mu(x,u,I^\delta
    W\star \Lambda^N_s(x,u))\big)f(x,u) +\int_{\cal U} f(x,v)M(x,u,v)dv
    \bigg\}\Lambda^N_s(dx,du) ds \label{defmnf}
  \end{multline}
  is a c\`adl\`ag $L^2$ martingale starting from $0$ with predictable
  quadratic variation
  \begin{multline}
    \label{quadran}
    \langle
    Z^{N,f}\rangle_t={1\over N}\int_0^t\int_{\bar{\cal X}\times {\cal
        U}} \bigg\{2 m(x,u)|\nabla_xf|^2+
    \big(\lambda(x,u)+\mu(x,u,I^\delta W\star
  \Lambda^N_s(x,u))\big)f^2(x,u) \\ +\int_{\cal U} f^2(x,v)M(x,u,v)dv
    \bigg\}\Lambda^N_s(dx,du) ds
  \end{multline}
\end{lem}

We assume

\noindent{\bf Assumption (H1):} \\
{\it 1) The initial measures $\Lambda^N_0$ converge in law and for
the weak
  topology on $M_F(\bar{\cal X}\times {\cal U})$ to some deterministic
  finite measure $\xi_0 \in M_F(\bar{\cal X}\times {\cal
    U})$, and $\sup_N E(\langle \Lambda_0^N,1\rangle^3)<+\infty$.}\\
{\it 2) All the parameters of the model are assumed to be
continuous,
  either on $\bar{\cal X}\times {\cal U}$, or on $\bar{\cal X}\times
  {\cal
    U}\times \rit$.}\\
{\it 3) There exists a constant $k_\mu$ such that
  \begin{equation}
    \label{death1} \forall x\in{\cal X},\:u\in{\cal U},\:r_1,
    r_2 \in \rit,  \quad | \mu(x,u,r_1)-\mu(x,u,r_2)|\leq
    k_\mu|r_1-r_2|.
  \end{equation}}

By the law of large numbers, Assumption (H1-1) is for example
satisfied for $\Lambda^N_0={1\over
N}\sum_{i=1}^N\delta_{(X^i_0,U^i_0)}$, with independent random
variables $\ (X^i_0,U^i_0)_{\{1\leq i\leq
  N\}}\ $ distributed following the law $\xi_0$ with finite 3rd-order moment.

Let us recall that the parameters of diffusion, birth and mutation
associated with $\Lambda^N$ stay unchanged, whereas the parameter
of selection $\mu_N$ is defined by (\ref{sn}).

\begin{thm}
  \label{conv.deter} Assume Hypotheses (H) and (H1), and consider the
  sequence of processes $\Lambda^N$ defined by (\ref{mvp}). Then for all
  $T>0$, the sequence $(\Lambda^N)$ converges in law, in
  $\dit([0,T],M_F(\bar{\cal X}\times~{\cal U}))$, to a deterministic
  continuous function $\: \xi^\delta\: $ belonging to
  $C([0,T],M_F(\bar{\cal X}\times {\cal U}))$.

  This measure-valued function $\ \xi^\delta\ $ is the unique weak
  solution satisfying $\ \sup_{t\in [0,T]} \langle\xi^\delta_t,1
  \rangle<+\infty\ $ of the following nonlinear integro-differential
  equation. For all function $f \in C^{2,b}_0$,
  \begin{align}
    \langle\xi^\delta_t,f\rangle & =\left<\xi_0,f\right> + \int_0^t
    \int_{\bar{\cal X}\times {\cal U}}\bigg\{m(x,u) \Delta_xf(x,u)+
    b(x,u).\nabla_x f(x,u)\notag\\ &+\big(\lambda(x,u) -\mu(x,u,I^\delta
    W\star \xi^\delta_s(x,u))\big)f(x,u) +\int_{\cal U}
    f(x,v)M(x,u,v)dv\bigg\}\ \xi^\delta_s(dx,du)  ds \label{dpde}
  \end{align}
\end{thm}

\begin{rema}
  \label{mass}
  Applying (\ref{dpde}) to the constant function equal to $1$, the
  positivity of $\mu$ and Hypotheses (H) gives
  $\langle\xi_t^\delta,1\rangle\ \leq \
  \langle\xi_0,1\rangle+C\int_0^t\langle\xi_s^\delta,1\rangle ds$.  We
  conclude by Gronwall's lemma that any solution $\xi^\delta$ of
  (\ref{dpde}) is bounded on every finite time interval $[0,T]$:
  $$
  \sup_{t\in [0,T]}\langle\xi^\delta_t,1\rangle\ \leq \ \langle\xi_0,1\rangle
  e^{C T}.
  $$
\end{rema}

%Before proving Theorem \ref{conv.deter}, we prove that each
%solution of (\ref{dpde}) is solution of an evolution equation.

As a first step in the proof of Theorem 4.2, we now give a mild
formulation for solutions of (\ref{dpde}). To this aim, and for
each fixed trait $u\in {\cal U}$, we denote by $\: P^{u}\: $ the
semigroup of the diffusion process normally reflected at the
boundary of ${\cal X}$, with diffusion matrix
$m(\cdot,u)\mbox{Id}$ and drift coefficient $b(\cdot,u)$.

\begin{lem}
  \label{evol} Let us consider a solution $\xi^{\delta}$ of
  (\ref{dpde}). Then, for each measurable and bounded function
  $\varphi$ defined on ${\bar{\cal X}}\times {\cal U}$,
  \begin{multline}
    \label{evolution}
    \langle \xi^\delta_t,\varphi\rangle =\langle
    \xi_0,P^{u}_{t}\varphi\rangle +\int_0^t\int_{\bar{\cal X}\times {\cal U}}
    \bigg\{\big(\lambda(x,u)-\mu(x,u,I^\delta
    W\star\xi^\delta_s(x,u))\big)P^{u}_{t-s}\varphi(x,u)\\ +\int_{\cal
      U} P^{v}_{t-s}\varphi(x,v)M(x,u,v)dv \bigg\}\
    \xi^\delta_s(dx,du)ds.
  \end{multline}
\end{lem}

\begin{proof}
  We may classically derive from (\ref{dpde}) a space-time weak
  equation for measurable functions $\psi_s(x,u)=\psi(s,x,u)$ which
  are of class $C^{1,2}$ on $[0,t] \times \bar{\cal X}$, measurable
  and bounded in $u$ and such that $\partial_n \psi=0$ on $[0,t]\times
  \partial {\cal X}\times {\cal U}$, given by
  \begin{align}
    \langle \xi^\delta_t,\psi_t\rangle & =
    \langle \xi_0,\psi_0\rangle+\int_0^t \int_{\bar{\cal X}\times {\cal U}}
    \bigg\{\partial_s \psi_s(x,u)+m(x,u)\Delta_x
    \psi_s(x,u)+b(x,u).\nabla_x\psi_s(x,u)\notag\\
    &+\big(\lambda(x,u)-\mu(x,u,I^\delta W\star
    \xi^\delta_s(x,u))\big)\psi_s(x,u) +\int_{\cal U}
    \psi_s(x,v)M(x,u,v)dv\bigg\}\ \xi^{\delta}_s(dx,du)ds
    \label{solmesd3}
  \end{align}
  Let us now consider a continuous function $\varphi$ on $\bar{\cal
    X}\times {\cal U}$ and fix a time $t\in[0,T]$. Let us define for
  $(s,x,u)\in [0,t]\times \bar{\cal X}\times {\cal U}$,
  $$\psi_s(x,u)=P^{u}_{t-s}\varphi(x,u).$$
  Then $\psi$ is solution of the  boundary value problem
  \begin{align*}
    &\partial_s \psi_s(x,u)+m(x,u)\Delta_x
    \psi_s(x,u)+b(x,u).\nabla_x\psi_s(x,u)=0\quad \hbox{ on
    }[0,T]\times {\cal X}\times {\cal U} \\ & \partial_n \psi_s(x,u)=0
    \quad \hbox{ on }[0,T]\times \partial{\cal X}\times {\cal U} \\
    &\psi_t(x,u)=\varphi(x,u) \quad \hbox{ on }\bar{\cal X} \times
    {\cal U}.
  \end{align*}
  Equation (\ref{solmesd3}) applied to this function $\psi$ yields the
  evolution equation
  \begin{multline}
    \langle
    \xi^\delta_t,\varphi\rangle =\langle \xi_0,P^{u}_{t}\varphi\rangle
    +\int_0^t\int_{\bar{\cal X}\times V}
    \bigg\{\big(\lambda(x,u)-\mu(x,u,I^\delta
    W\star\xi^\delta_s(x,u))\big)P^{u}_{t-s}\varphi(x,u)\\ +\int_{\cal
      U} P^{v}_{t-s}\varphi(x,v)M(x,u,v)dv \bigg\}\
    \xi^\delta_s(dx,du)ds.  \label{evolution2}
  \end{multline}
  Equation (\ref{evolution2}) is true for each continuous (and then
  bounded) function $\varphi$, and characterizes the finite measure
  $\xi^\delta$. Lemma \ref{evol} is proved.
\end{proof}

\begin{proof} (of Theorem \ref{conv.deter}).  Let us fix
  $T>0$.\\
  Let us firstly prove the uniqueness of solutions $\xi$ of
  (\ref{dpde}). Using Remark \ref{mass} and Lemma \ref{evol}, we prove
  the uniqueness of bounded solutions of (\ref{evolution}).  Let us
  consider two such solutions $(\xi_t)_{t\geq 0}$ and $(\bar
  \xi_t)_{t\geq 0}$ and compute the quantity $|\langle\xi_t-\bar
  \xi_t,\varphi\rangle|$, for each measurable and bounded function
  $\varphi$ such that $\|\varphi\|_\infty\leq 1$.

  Using (\ref{evolution}), we obtain for $t\leq T$
  \begin{align*}
    &\vert \left<\xi_t-\bar \xi_t,\varphi\right> \vert \leq \int_0^t
    \bigg|\int_{\bar{\cal X}\times {\cal U}}
    \bigg\{\big(\lambda(x,u)-\mu(x,u,I^\delta W\star
    \xi_s(x,u))\big)P^{u}_{t-s}\varphi(x,u) \\ & 
    \phantom{\vert \left<\xi_t-\bar \xi_t,\varphi\right> \vert \leq \int_0^t
    \bigg|\int_{\bar{\cal X}\times {\cal U}}
    \bigg\{ }
  +\int_{\cal U}
    P^{v}_{t-s}\varphi(x,v)M(x,u,v)dv\bigg\}\left(\xi_s(dx,du) - \bar
      \xi_s(dx,du)\right)\bigg| ds \\ &+  \int_0^t \int_{\bar{\cal
        X}\times {\cal U}} \left\vert \big(\mu(x,u,I^\delta W\star \bar\xi_s(x,u))
      -\mu(x,u,I^\delta W\star \xi_s(x,u))\big)P^{u}_{t-s}\varphi(x,u)
    \right\vert\bar \xi_s(dx,du) ds
  \end{align*}
  Now, using Hypotheses (H), applying Remark \ref{mass} to $\bar{\xi}$
  and since $\|\varphi \|_\infty \leq 1$, there exists a positive
  constant $C_1$ such that for all $(x,u)\in \bar{\cal X}\times {\cal
    U}$ and all $0<s\leq t\leq T$,
  \begin{align*}
   |\lambda(x,u)P^{u}_{t-s}\varphi(x,u)+\int_{\cal U}
    P^{v}_{t-s}\varphi(x,v)M(x,u,v)dv| & \ \leq\ C_1,\\
   |\mu(x,u,I^\delta W\star
    \bar\xi_s(x,u))P^{u}_{t-s}\varphi(x,u)|\leq \mu_0(1+\|I^\delta
    W\|_\infty\langle\bar\xi_s,1\rangle{}) & \ \leq \ C_1
  \end{align*}
  while thanks to (H1-2),
  $$|\mu(x,u,I^\delta
  W\star \bar\xi_s(x,u)) -\mu(x,u,I^\delta W\star \xi_s(x,u)) |\leq
  k_\mu\|I^\delta W\|_\infty \sup_{\|\varphi\|_\infty\leq
    1}|\langle\xi_s - \bar \xi_s,\varphi\rangle|,$$
  and then
  \begin{align*}
   & \left|\int_{\bar{\cal X}\times {\cal U}} \left(\mu(x,u,I^\delta
        W\star \bar\xi_s(x,u)) -  \mu(x,u,I^\delta W\star
        \xi_s(x,u))\right)P^{u}_{t-s}\varphi(x,u)\bar
      \xi_s(dx,du)\right| \\ & \hskip 7cm \leq C_2\sup_{\|\varphi\|_\infty\leq
      1}|\langle{}\xi_s - \bar \xi_s,\varphi\rangle{}|
  \end{align*}
  where $C_2$ is a positive constant. We deduce that there exists
  $C>0$ such that
  \begin{equation*}
    \vert \left<\xi_t-\bar \xi_t,\varphi\right> \vert \leq C \int_0^t
    \sup_{\|\varphi\|_\infty\leq 1}|\langle\xi_s - \bar \xi_s,\varphi\rangle|ds
  \end{equation*}
  and by Gronwall's lemma, we conclude that for all $t\leq T$,
  $
  \sup_{\|\varphi\|_\infty\leq 1}|\langle\xi_t - \bar \xi_t,\varphi\rangle|=0.
  $
  Thus, for all $t\leq T$, $ \xi_t = \bar\xi_t$ and uniqueness holds.

  Let us next prove that for all $T>0$,
  \begin{equation}
    \label{mom3}
    \sup_{N\in \nit^*} E\left(\sup_{[0,T]} \left< \Lambda^N_t,1\right>^3\right) <+\infty
  \end{equation}
  Introducing $S_k^N=\inf\{t\geq 0, \left< \Lambda^N_t,1\right>\geq k\}$
  for $k\in\nit^*$, a simple computation using the specific form of
  $L_{N,1}F_f$ and $L_{N,2}F_f$ with $f=1$ and $F(y)=y^3$ and
  dropping the negative death term yields
  \begin{equation*}
    E\left(\sup_{s\leq
        t\wedge S_k^N}\left< \Lambda^N_s,1\right>^3\right)\leq E(\left<
      \Lambda^N_0,1\right>^3)+CE\left(\int_0^{t\wedge S^N_k}(\langle
      \Lambda^N_s,1\rangle+\left< \Lambda^N_s,1\right>^3)ds \right)
  \end{equation*}
  where $C$ is a positive constant independent of $k$ and $N$. Then Assumption
  (H1-1) and Gronwall's lemma imply that there exists a constant $C_T$
  independent of $k$ and $N$ such that $ E\left(\sup_{s\leq T\wedge
      S_k^N}\left< \Lambda^N_s,1\right>^3\right)\leq C_T.$ We deduce that
  the sequence $(S_k^N)_k$ tends a.s.\ to infinity and finally obtain
  (\ref{mom3}) by Fatou's lemma.

  Using Remark \ref{densite}, and following Roelly \cite{Roelly:86},
  one observes that the sequence of laws $Q^N$ of $\Lambda^N$ is uniformly
  tight in ${\cal P}(\dit([0,T],M_F(\bar{\cal X}\times {\cal U})))$,
  where $M_F$ is endowed with the vague topology, as soon as for any
  function $f\in C^{2,0}_0$, the sequence of the laws of the processes
  $\left<\Lambda^N,f\right>$ is tight in ${\cal P}(\dit([0,T], \rit))$.
  Using Aldous' \cite{Aldous:78} and Rebolledo's \cite{Joffe:86}
  criteria, this tightness follows from
  \begin{equation}
    \label{tight} \sup_{N \in \nit^*}
    E(\sup_{[0,T]} \vert \left< \Lambda^N_s,f\right> \vert) <\infty,
  \end{equation}
  and from the tightness of the laws of $(\langle Z^{N,f}\rangle)$ and
  of the drift part of the semimartingales $\left<\Lambda^N,f\right>$. \\
  Clearly, since $f$ is bounded, (\ref{tight}) is a consequence of
  (\ref{mom3}). Let us now consider stopping times $(S,S')$ satisfying
  a.s. $0 \leq S \leq S' \leq S+\delta\leq T$. Thanks to Doob's
  inequality, Lemma \ref{martn}, and (\ref{mom3}), we get
  \begin{equation*}
    E\left(\langle Z^{N,f}\rangle_{S'}-
      \langle Z^{N,f}\rangle_S\right)\leq E\left( C \int_S^{S+\delta}
      \left( \left< \Lambda^N_s,1\right>+ \left< \Lambda^N_s,1\right>^2
      \right)ds\right) \leq C \delta.
  \end{equation*}
  Similar arguments prove that the expectation of the finite variation
  part of $ \left< \Lambda^N_{S'},f\right> - \left< \Lambda^N_S,f\right>$ is
  bounded by $C \delta$.
  Finally it turns out that the sequence
  $(Q^N)_N$ is uniformly tight.

  Let us now denote by $Q$ the limiting law in ${\cal
    P}(\dit([0,T],M_F(\bar{\cal X}\times {\cal U})))$ of a subsequence
  of $Q^N$, still denoted by $Q^N$ for simplicity. By construction,
  almost surely,
  \begin{equation*}
    \sup_{t\in [0,T]}  \sup_{\vert \vert f \vert \vert_\infty\leq 1}
    \vert \langle \Lambda^N_s,f \rangle - \langle \Lambda^N_{s-},f\rangle \vert
    \leq 1/N.
  \end{equation*}
  We deduce immediately that each process $\Lambda$ with law $Q$ is
  a.s.\ strongly continuous.
  Let us finally prove that it is the unique solution of
  (\ref{dpde}) .

  For $t\leq T$, $f\in C_0^{2,0}$ and $\nu \in
  \dit([0,T],M_F(\bar{\cal X}\times {\cal U}))$, let us define
  \begin{align*}
    \Psi_{t,f}(\nu)&= \left<\nu_t,f\right> - \left<\nu_0,f\right> -
    \int_0^t \int_{\bar{\cal X}\times {\cal U}} \bigg\{m(x,u)
    \Delta_xf(x,u)+b(x,u).\nabla_xf(x,u) \\ &+ \big(\lambda(x,u)
    -\mu(x,u,I^\delta W\star \nu_s(x,u))\big)f(x,u) +\int_{\cal U}
    f(x,v)M(x,u,v)dv\bigg\}\nu_s(dx) ds.
  \end{align*}
  We want to show that for any $t\leq T$,
  \begin{equation}\label{wwhtp}
    E \left( |\Psi_{t,f}(\Lambda) | \right)=0,
  \end{equation}
  knowing from Lemma \ref{martn} that
  \begin{equation}\label{wwk}
    Z^{N,f}_t = \Psi_{t,f}(\Lambda^N).
  \end{equation}
  A fair computation using Lemma \ref{martn}, Hypotheses (H) and (H1),
  and (\ref{mom3}) shows that
  \begin{equation}\label{cqmke}
    E \left( \vert Z^{N,f}_t \vert^2 \right) = E \left( \langle
      Z^{N,f}\rangle_t \right) \leq \frac{C_f}{N}E\left( \int_0^t
      \left\{1+\left< \Lambda^N_s,1\right>^2\right\} ds \right) \leq
    \frac{C_{f,t}}{N}
  \end{equation}
  which goes to $0$ as $N$ tends to infinity. On another hand, since
  $\Lambda$ is a.s.\ strongly continuous, since $f\in C_0^{2,0}$ and thanks
  to the assumption (H), the function $\Psi_{t,f}$ is a.s.\ continuous
  at $\Lambda$. Furthermore, for any $ \nu \in \dit([0,T],M_F(\bar{\cal
    X}\times {\cal U}))$,
  \begin{equation*}
    \vert \Psi_{t,f}(\nu)\vert \leq C_{t,f} \sup_{[0,t]}
    \left(1+\left< \nu_s,1\right>^2 \right),
  \end{equation*}
  and (\ref{mom3}) implies that the sequence $(\Psi_{t,f}(\Lambda^N))_N$ is
  uniformly integrable. Thus
  \begin{equation}\label{cqv1}
    \lim_N E \left( |\Psi_{t,f}(\Lambda^N)|\right)= E \left( |\Psi_{t,f}(\Lambda)|
    \right)
  \end{equation}
  Combining (H1-1), (\ref{wwk}), (\ref{cqmke}) and (\ref{cqv1}), we
  conclude that (\ref{wwhtp}) holds and that (\ref{dpde}) is satisfied
  for any $f\in C^{2,0}_0$.

  Then $\Lambda$ is uniquely identified to $\xi^{\delta}$, and
 the sequence $(\Lambda^N)$ converges to
  $\xi^{\delta}$ in $\dit([0,T], M_F(\bar{\cal X}\times~{\cal U}))$,
  where $M_F(\bar{\cal X}\times {\cal U})$ is endowed with the vague
  topology. To extend this result to the weak topology, we use a
  criterion proved in \cite{Meleard:93}. Since the limiting process is
  continuous, it suffices to prove that the sequence $(\langle \Lambda^N,
  1\rangle)_N$ converges in law to $\langle \xi^{\delta}, 1\rangle$ in
  $\dit([0,T], \rit)$. We may apply what has been done above with
  $f\equiv 1$. Theorem \ref{conv.deter} is proved.
\end{proof}

In the next section, we will be interested  in the limit of small
spatial interaction range $\delta$. An intermediate result
consists in proving the existence of a density for each measure
$\xi^{\delta}_t$, $t\geq 0$. We make the
additional\\
{\bf Assumption (H2):}\\
{\it 1) The diffusion coefficient $m(x,u)$ is of class $C^2$ in
$x$
  and the second derivative of $m$ (in $x$) is $\alpha$-H\"olderian,
  uniformly in $u$, with $\alpha>0$. Moreover, $m$ is assumed to be
  positive. Hence, since $\bar{\cal X}\times {\cal U}$ is a compact
  set, there exists $m_*>0$ such that for all $(x,u)\in \bar{\cal
    X}\times U$,
  $$ m(x,u)\geq m_*>0.$$
  2) The drift coefficient $b(x,u)$ is of class $C^1$ in $x$ and the
  derivative of $b$ (in $x$) is $\alpha$-H\"olderian, uniformly in
  $u$, with $\alpha>0$.}

Assumptions (H) and (H2) and the smoothness of $\partial {\cal X}$
allow us to adapt Sato-Ueno \cite{Sato:65} (Theorem 2.1 and
Appendix) to obtain the following lemma.

\begin{lem}
  \label{semigroup} There exists a unique function $p_t(x,u,y)$
  defined on $\rit_+\times {\bar{\cal X}}\times {\cal U}\times
  {\bar{\cal X}}$, continuous in $(t,x,y)$ and which is a density
  function in $y\in \bar{\cal X}$ such that for each continuous
  function $\varphi$ defined on $\bar{\cal X}\times U$, each $(x,u)\in
  \bar{\cal X}\times {\cal U}$,
  \begin{equation} \label{sg}
    P^u_t\varphi(x,u)=\int_{{\bar{\cal X}}}p_t(x,u,y)\varphi(y,u)dy
  \end{equation}
\end{lem}
Let us now prove the propagation in time of the absolute
continuity property of the measure-valued solution $\xi^{\delta}$.

\begin{thm}
  \label{density} Assume (H), (H1) and (H2) and that
  $\xi_0(dx,du)=g_0(x,u)dxdu$. Then for each time $t$, the measure
  $\xi^\delta_t$ has a density $g^{\delta}\in L^{\infty}([0,T],L^1)$
  with respect to the Lebesgue measure on $\bar{\cal X}\times {\cal
    U}$. Moreover, for each $t$ and $u$, the function
  $g^{\delta}_t(.,u)$ is continuous on ${\bar{\cal X}}$.
\end{thm}

\begin{proof} Let us come back to the equation (\ref{evolution})
  satisfied by $\xi^\delta$.

%We construct recursively an implicit sequence $(g_n)_n$ of
%positive functions belonging to $L^{\infty}_T(L^1(\bar{\cal
%X}\times U))$ and which will converge to a
%$L^{\infty}_T(L^1(\bar{\cal X}\times U)$-function $g$ being a weak
%solution of (\ref{evolution}). By uniqueness, the relut follows.

  Using basic results on linear parabolic equations, we construct by
  induction a sequence of functions $(g_n)_n$ satisfying in a weak sense
  \begin{align}
    &\partial_t
    g^{n+1}_t(x,u)=\Delta_x(m(x,u)g^{n+1}_t(x,u))-\nabla_x(b(x,u)g^{n+1}_t(x,u))
    \notag\\
    & \qquad+\lambda(x,u)g^n_t(x,u) +\int_{\cal U} g_t^n(x,v) M(x,u,v)dv
    -\mu(x,u,I^\delta W\star g^n_t(x,u)) g^{n+1}_t(x,u)\notag\\
    &g^{n+1}_0(x,u) =g_0(x,u)\notag\\ &\nabla_x g^{n+1}(t,x,u).n(x)=0
    \quad \forall (t,x,u)\in \rit_+\times
    \partial {\cal X}\times {\cal U}.\label{eqn.gn}
  \end{align}
  Thanks to the nonnegativity of $g_0$, $\mu$, $\lambda$ and
  $M$, and applying the maximum principle, we can show that the
  functions $g_n$ are nonnegative (see \cite{Desville:04}). By symmetry of $M$, Equation
  (\ref{eqn.gn}) (understood in the weak sense) means that for all
  $C^{2,b}_0$-function $\varphi$ from $\bar{\cal X}\times {\cal U}$
  into $\rit$,
  \begin{align}
    \langle g^{n+1}_t,\varphi\rangle = \left<g_0,\varphi\right> &+
    \int_0^t \int_{\bar{\cal X}\times {\cal U}}\bigg\{\bigg(m(x,u)
    \Delta_x\varphi(x,u)+ b(x,u).\nabla_x
    \varphi(x,u)\bigg)g^{n+1}_s(x,u)\notag\\
    &+\bigg(\lambda(x,u)\varphi(x,u) +\int_{\cal U}
    \varphi(x,v)M(x,u,v)dv\bigg)g^n_s(x,u)\notag\\ &-\mu(x,u,I^\delta
    W\star g^n_s(x,u))\varphi(x,u)g^{n+1}_s(x,u)\bigg\} dx\: du\: ds.\label{dpgn}
  \end{align}

  The associated mild equation writes as before: for each continuous
  function $\varphi$,
  \begin{align}
    \langle g^{n+1}_t,\varphi\rangle & =\int_{\bar{\cal X}\times
      {\cal U}}\bigg(\int_{\bar{\cal X}}
    p_t(x,u,y)\varphi(y,u)dy\bigg)g_0(x,u)dxdu\notag\\ &+\int_0^t
    \int_{\bar{\cal X}\times {\cal U}}
    \bigg\{\bigg[\lambda(x,u)\bigg(\int_{\bar{\cal X}}
    p_{t-s}(x,u,y)\varphi(y,u)dy\bigg) \notag\\ &+\int_{\cal U}
    \bigg(\int_{\bar{\cal X}}
    p_{t-s}(x,v,y)\varphi(y,v)dy\bigg)M(x,u,v)dv\bigg]g^n_s(x,u)\notag\\
    &-\mu(x,u,I^\delta
    W\star g^n_s(x,u))\bigg(\int_{\bar{\cal X}}
    p_{t-s}(x,u,y)\varphi(y,u)dy\bigg)g^{n+1}_s(x,u)\bigg\} dx du
    ds.\label{evolgn}
  \end{align}

  Hypotheses on the coefficients allow us to apply Fubini's theorem
  and to obtain that for each $(y,u)\in \bar{\cal X}\times {\cal U}$,
  \begin{align}
    g^{n+1}_t(y,u)
    &=\int_{\bar{\cal X}} p_t(x,u,y)g_0(x,u)dx\notag\\
    &+\int_0^t\int_{\bar{\cal X}}
    \bigg \{(\lambda(x,u)p_{t-s}(x,u,y) g^n_s(x,u) +\int_{\cal U}
    p_{t-s}(x,u,y) g^n_s(x,v)M(x,u,v)dv\notag\\ & 
    \phantom{+\int_0^t\int_{\bar{\cal X}} \bigg \{ }
    -\mu(x,u,I^\delta
    W\star g^n_s(x,u)) p_{t-s}(x,u,y)g^{n+1}_s(x,u)\bigg\}dxds.
    \label{induction}
  \end{align}
  Then, thanks to the nonnegativity of $g^{n+1}$, we get
  \begin{multline}
    0 \leq g^{n+1}_t(y,u)\leq \int_{\bar{\cal X}}
    p_t(x,u,y)g_0(x,u)dx +\int_0^t\int_{\bar{\cal X}} \bigg
    \{(\lambda(x,u)p_{t-s}(x,u,y) g^n_s(x,u) \\ +\int_{\cal U}
    p_{t-s}(x,u,y) g^n_s(x,v)M(x,u,v)dv\bigg\}dxds, \label{pos}
  \end{multline}
  and deduce easily, integrating over $y\in\bar{\cal X}$, using
  Fubini's Theorem, the symmetry of $M$ and Gronwall's Lemma that
  there exists a constant $C$ independent of $\delta$ such that
  \begin{equation} \label{major.gn} \sup_{n\in \nit}\sup_{t\leq
      T}\|g^n_t\|_1 \leq \|g_0\|_1e^{C T}.
  \end{equation}

  Let us now show the convergence of the sequence $g^n$ in
  $L^{\infty}([0,T],L^1)$ to a function $g^{\delta}$. A
  straightforward computation using (\ref{induction}), Hypotheses (H),
  (H1) and (H2), and similar arguments as above yields
  $$
  \sup_{s\leq t}\|g^{n+1}_s-g^n_s\|_1 \leq
  C\int_0^t\bigg(\sup_{u\leq s}\|g^{n+1}_u-g^n_u\|_1+\sup_{u\leq
    s}\|g^{n}_u-g^{n-1}_u\|_1\bigg) ds
  $$
  where $C$ is a positive constant. Thanks to Gronwall's Lemma, we
  deduce that for each $T>0$, each $t\leq T$ and each $n$,
  $
  \sup_{s\leq t}\|g^{n+1}_s-g^n_s\|_1  \leq C\int_0^t\sup_{u\leq
    s}\|g^{n}_u-g^{n-1}_u\|_1 ds$.\\
  Picard's Lemma allows us to conclude that for any $T>0$,
  $\ \sum_{n\in\nit}\sup_{t\in[0,T]}\|g^{n+1}_t-g^n_t\|_1<+\infty,\ $
  and the sequence $(g^n)_n$ converges in $L^{\infty}([0,T],L^1)$ to a
  function $g^{\delta}$.  We deduce from (\ref{major.gn}) that
  \begin{equation}
    \label{major.gd} \sup_{\delta>0}\sup_{t\leq T}\|g^{\delta}_t\|_1
    \leq \|g_0\|_1e^{C T}.
  \end{equation}

  Moreover, the function $g^{\delta}$ is solution of (\ref{dpde}), and
  thus, the uniqueness result proved in Theorem \ref{conv.deter}
  implies that $\xi^{\delta}(dx,du)=g^{\delta}(x,u)dxdu$. Then, the
  measure $\xi^{\delta}$ is absolutely continuous with respect to the
  Lebesgue measure, and the density $g^\delta$ is weak solution of the
  nonlocal nonlinear partial differential equation
  \begin{align}
    &\partial_t
    g^\delta=\Delta_x(m(x,u)g^\delta_t(x,u))-\nabla_x(b(x,u)g^\delta_t(x,u))
    \notag\\ &\qquad+\big(\lambda(x,u)-\mu(x,u,I^\delta W\star
    g^{\delta}_t(x,u))\big)g^{\delta}_t(x,u)+\int_{\cal U}
    g_t^{\delta}(x,v)M(x,u,v)dv\ ;\notag\\ &g^{\delta}_0(x,u)
    =g_0(x,u)\ ;\notag\\ &\nabla_x g^\delta(t,x,u).n(x)=0  \quad \forall
    (t,x,u)\in \rit_+\times
    \partial {\cal X}\times {\cal U}.\label{w.gd}
  \end{align}
  Lemma \ref{evol} implies that $g^{\delta}$ is also solution of the
  mild equation
  \begin{align}
    g_t^\delta(y,u) & =\int_{{\cal X}}
    p_t(x,u,y)g_0(x,u)dx\notag\\ &+\int_0^t \int_{{\cal X}}
    \bigg\{\big(\lambda(x,u)- \mu(x,u,I^\delta
    W\star g^\delta_s(x,u))\big)p_{t-s}(x,u,y)g^\delta_s(x,u)\notag\\
    & \phantom{+\int_0^t \int_{{\cal X}} \bigg\{ }
    + \int_{\cal U}
    p_{t-s}(x,u,y)g^\delta_s(x,v)M(x,u,v)dv \bigg\}dxds.
    \label{evol.gd}
  \end{align}
  Using (\ref{evol.gd}), the continuity of $y\mapsto g^\delta_t(y,u)$
  follows immediately from the continuity of $(t,x,y)\mapsto
  p_t(x,u,y)$, the nonnegativity and boundedness of $g^\delta$ and the
  boundedness of birth parameters.
\end{proof}

\setcounter{equation}{0}
\section{Convergence of the number density when the interaction range
  decreases}
\label{sec:delta=0}

Our aim in this section is to prove that under suitable assumptions,
the sequence $(g^\delta)$ converges, when $\delta$ tends to $0$, to a
function $g\in L^{\infty}([0,T],L^1)$ with initial condition $g_0$,
which is weak solution of the locally nonlinear partial differential
equation
\begin{align}
  &\partial_tg_t(x,u)
  =\Delta_x(m(x,u)g_t(x,u))-\nabla_x(b(x,u)g_t(x,u))\notag\\ &
  \phantom{\partial_tg_t(x,u)}
  +\big(\lambda(x,u)-\mu(x,u,\rho_g(t,x,u))\big)g_t(x,u)+\int_{\cal U}
  g_t(x,v)M(x,u,v)dv\ ; \notag\\
  &\nabla_x g(t,x,u).n(x)=0 \quad \forall (t,x,u)\in \rit_+\times
  \partial {\cal X}\times {\cal U} \label{equa}
\end{align}
where $\rho_g$ describes the (local) interaction in $x$, defined for
$(x,u)\in \bar{\cal X}\times {\cal U}$ by
$$
\rho_g(t,x,u)=\int_{\cal U} W(u-v)g_t(x,v)dv.
$$
In order to control the terms $I^{\delta}W*g^{\delta}$ uniformly
in $\delta$ in the nonlinear term of (\ref{w.gd}), we need
$L^{\infty}$-estimates on $g^{\delta}$ and we make the following
initial
data assumption:\\\\
{\bf (H3)\ }{\it The initial density $g_0(x,u)$ is bounded on ${\cal X}\times
  {\cal U}$.}

\begin{prop}
  \label{unid}
  Assume (H), (H1), (H2), (H3). Then there exists a positive constant
  $C_T$, such that
  \begin{equation} \label{unifd}
    \sup_{\delta>0}
    \sup_{t\in[0,T]}\|g^\delta_t\|_\infty \leq C_T\|g_0\|_{\infty}.
  \end{equation}
\end{prop}

\begin{proof}
  Let us again consider the sequence $(g^n)_n$ approximating
  $g^{\delta}$ introduced in the proof of Theorem \ref{density}. The
  maximum principle implies that
  $$\sup_n\sup_{t\leq T}\|g^n_t\|_{\infty}\leq C\|g_0\|_{\infty},$$
  where $C>0$ is a constant only depending on $T$, $\lambda^*$ and
  $M^*$ (and independent of $\delta$). This property propagates
  taking the limit in $n$, and (\ref{unifd}) is proved. (For details
  on the maximum principle, see \cite{Desville:04}.)
\end{proof}

Let us now prove the following convergence theorem:
\begin{thm}
  \label{convergence}
  Assume hypotheses (H), (H1), (H2), (H3).  Assume that the measure
  $I^\delta(y)dy$ weakly converges to the Dirac measure $\delta_0$ as
  $\delta$ tends to $0$. (To fix ideas we may assume that
  $I^{\delta}(x)=C_{\delta}{\bf 1}_{\{|x|\leq \delta\}}$.) Then the
  sequence $(g^\delta)_{\delta>0}$ converges in
  $L^{\infty}([0,T],L^1)$ as $\delta$ tends to $0$, to the unique
  function $g \in L^{\infty}([0,T], L^1\cap L^{\infty}(\bar{\cal
    X}\times {\cal U}))$ satisfying for each $y,u \in \bar{\cal X}\times
  {\cal U}$ the evolution equation
  \begin{align}
    g_t(y,u) &=\int_{{\cal X}} p_t(x,u,y)g_0(x,u)dx\notag\\
    &+\int_0^t \int_{{\cal X}}
    \bigg\{\big(\lambda(x,u) - \mu(x,u,\rho_g(s,x,u))\big)
    p_{t-s}(x,u,y)g_s(x,u)\notag\\ & \qquad\qquad\ + \int_{\cal U}
    p_{t-s}(x,u,y)g_s(x,v)M(x,u,v)dv
    \bigg\}dxds. \label{evol.g}
  \end{align}
  Moreover, for each $t$ and $u$, the function $g_t(.,u)$ is
  continuous on ${\cal X}$.
\end{thm}

% Observe that, since ${\cal X}$ is bounded, the natural renormalization
% of the interaction kernel as $I^\delta(y)=I(y/\delta)/\delta^d$ yields
% an spatial interaction kernel that satifies the assumption on the
% diameter of the support of $I^\delta$.

\begin{proof}
  One can easily prove the existence and uniqueness of the integrable
  and bounded function $g$ solution of (\ref{evol.g}) by adapting the
  proofs of Theorem \ref{density} and Proposition \ref{unid},
  replacing $\mu(x,u,I^\delta W\star~g)$ by $\mu(x,u,\rho_g)$.
  The continuity of $y\to g_t(y,u)$ is obtained as in the proof of Theorem
  \ref{density}, and we can show as in the proof of Proposition
  \ref{unid} that
  \begin{equation}
    \label{eq:borne-g}
    \sup_{t\in[0,T]}\|g_t\|_{\infty}\leq C_T\|g_0\|_{\infty}.
  \end{equation}

  Let us write
  \begin{align}
    g^\delta_t(y,u) & -g_t(y,u)=\int_0^t
    \int_{{\cal X}}
    \bigg\{(\lambda(x,u)p_{t-s}(x,u,y)\big(g^\delta_s(x,u)-g_s(x,u)\big)
    \notag \\
    &+ \int_{\cal U}
    p_{t-s}(x,u,y)\big(g^\delta_s(x,v)-g_s(x,v)\big)M(x,u,v)dv \notag \\ &-
    \bigg[\mu(x,u,I^\delta W\star g^\delta_s(x,u))g^\delta_s(x,u)-
    \mu(x,u,\rho_g(s,x,u))g_s(x,u)\bigg]p_{t-s}(x,u,y)
    \bigg\}dxds \label{gdg}
  \end{align}
 Using (\ref{unifd}) and (\ref{eq:borne-g}), the unique term which deserves attention is the term
  $\mu(x,u,\rho_g(x,u))-\mu(x,u,I^\delta W\star g^\delta(x,u))$.
  By (\ref{death1}), we have
  \begin{align*}
    &\int_{\cal
      X}|\mu(x,u,\rho_g(t,x,u))-\mu(x,u,I^\delta W\star
    g^\delta_t(x,u))|dx \\ &\leq k_\mu\int_{\cal X}\left|\int_{\cal U}
      W(u-v)g_t(x,v)dv-\int_{{\cal X}\times {\cal
          U}}I^\delta(x-z)W(u-v)g^\delta_t(z,v)dzdv\right|dx \\ &\leq
    k_\mu\int_{\cal X}\bigg(\left|\int_{\cal U}
      W(u-v)g_t(x,v)dv-\int_{{\cal X}\times {\cal
          U}}I^\delta(x-z)W(u-v)g_t(z,v)dzdv\right| \\
    &\qquad+\int_{{\cal X}\times {\cal
        U}}I^\delta(x-z)W(u-v)\left|g_t(z,v)-g^\delta_t(z,v)\right|dzdv\bigg)dx
%     \\
%     &\leq k_\mu\bigg(\int_{\cal X}\int_{{\cal X}\times {\cal U}}
%     I^\delta(x-z)U(u-w)\left|g_t(x,w)-g_t(z,w)\right|dzdwdx\bigg) \\
%     &+\int_{\cal X}\int_{{\cal X}\times {\cal
%         U}}I^\delta(x-z)U(u-w)\left|g_t(z,w)-g^\delta_t(z,w)\right|dzdwdx\bigg)
  \end{align*}
  Let us fix our attention on the first term in the last right
  inequality, that we will call $A_{\delta}(t,u)$. Since
  $I^\delta(y)dy$ weakly converges to $\delta_0$, $\int_{{\cal
      X}}I^\delta(x-z)g_t(z,v)dz$ converges to $g_t(x,v)$ as $\delta$
  goes to 0. Because of~(\ref{eq:borne-g}), this convergence holds in
  a bounded pointwise sense with respect to $t\leq T$, $x\in{\cal X}$
  and $v\in{\cal U}$.
%  Since the range of
%   interaction is of order $\delta$ and using the uniform continuity of
%   $ g_t(.,u)$ on the compact set $\bar{\cal X}$, we get
%   $|g_t(x,w)-g_t(z,w)|\leq \eta(\delta,t,w)$ when $|x-z|\in {\rm
%     supp}(I^{\delta})$, where $\eta(\delta,t,w)$ tends to $0$ as
%   $\delta$ tends to $0$, in a bounded pointwise sense (using the
%   boundedness of $g$).
  Then Lebesgue's theorem implies that $A_{\delta,T}:=\int_{{\cal
      U}}\int_0^T A_{\delta}(t,u)dtdu$ tends to $0$ as $\delta$ tends
  to $0$.

  Now, integrating (\ref{gdg}) with respect to $dy\:du$, a
  straightforward computation yields
  \begin{equation*}
    \sup_{s\leq
      t}\|g^\delta_s-g_s\|_1\leq C_T A_{\delta,T} +
    C'_T\int_0^t \sup_{u\leq s}\|g^\delta_u-g_u\|_1ds.
  \end{equation*}
  We conclude using Gronwall's lemma.
\end{proof}

The zero interaction range equation (\ref{equa}) has been
numerically studied in Pr\'evost \cite{Prevost:04}.  A lot of
simulations based on finite element schemes are given, studying
the simultaneous effects of the diffusion, mutation and selection
on the invasion of the domain by the population. The simulations
show that the coefficient which seems to affect the most the
invasion aptitude is the mutation size coefficient. However, they
restrict to local interactions.

In the next section, we wish additionally to illustrate, by
simulations of the stochastic discrete model, the effect of the
spatial interaction range on the interplay between invasion and
evolution, and the emergence of spatial and phenotypic diversity
(clustering and polymorphism). Our simulations focus on the
qualitative differences between local and nonlocal interactions.

\setcounter{equation}{0}
\section{Simulations}
\label{sec:simul}

% The role of space in the coexitence of several types of individual and
% in the emergence and stability of polymorphism in a population, that
% may eventually lead to speciation, has long been recognized in the
% biological literature (Mayr~\cite{Mayr63}, Endler~\cite{Endler77},
% Durrett and Levin~\cite{DL94}, Dieckmann and Doebeli~\cite{DD99}). In
% particular, the role of space-related traits (such as the dispersal
% speed or the sensibility to heterogeneously distributed resources, see
% Bolker and Pacala~\cite{BP:99}) are often closely related to the
% appearance and stabilisation of polymorphism. Moreover, this
% polymorphism is often linked with some spatial clustering
% (organisation of the population as isolated ``patches'') of the
% population, linked to spatial heterogeneity (Flierl et
% al.~\cite{FGLO99}) or simply linked to the fact that reproductions are
% local (Young et al.~\cite{YRS01}). Here, we want to address more
% specifically, through some simulations, the role of the interaction
% range (parameter $\delta$ in our model), with respect to spatial
% structuration and polymorphism in the population. In particular, we
% want to show that in some situations, assuming as in (\ref{equa}) and
% \cite{Prevost:04} that spatial interactions are local may induce a
% loss in the detail of the population structure that can lead to a
% misunderstanding of the mechanisms at the origin of diversity.

We will give in this section simulations of several biologically
realistic examples, based on the algorithm of
Section~\ref{sec:mart}. The Euler scheme to simulate reflected
diffusions will be detailed in Section~\ref{sec:algo}, as well as
some simplifications in the algorithm of Section~\ref{sec:mart},
in the case of linear death rates.

Next, we will give simulations of three biologically relevant
examples.  First (Section~\ref{sec:ex1}), we show that, when
migrations and mutations are not too strong, a large interaction range
induces a spatial organization of the population as a finite set of
isolated clusters.  Conversely, for sufficiently small interaction
range, the clustering phenomenon is no more observed.  Second
(Section~\ref{sec:ex2}), we propose another example where a similar
phase transition occurs for spatial clustering and in which the
critical interaction range can be identified. In our last example
(Section~\ref{sec:invas}), we investigate a model describing the
invasion of a species with evolving dispersal speed (the trait is
proportional to the migration speed, as in \cite{Desville:04}).

\subsection{Euler scheme and algorithm for logistic interaction}
\label{sec:algo}

As mentionned in Section~\ref{sec:mart}, the reflected diffusion of
our particles can be simulated with an Euler scheme. We will assume in
this subsection and in the following examples that ${\cal
  X}=(\alpha,\beta)\subset\mathbb{R}$ and we will use the scheme of
L\'epingle~\cite{Le95} (see also~\cite{Gobet:01}). Fix
$x\in[\alpha,\beta]$ and $u\in{\cal U}$. On any time interval where
its trait is constant, an indivual at $(x,u)$ moves according to the
reflected diffusion
\begin{align}
  X_t & =x+\int_0^t
  \sqrt{2m(X_{s},u)}dB_s
  +\int_0^t b(X_s,u)ds -k_t\ , \notag\\
  |k|_t & =\int_0^t \1_{\{X_s\in\{\alpha,\beta\}\}}d|k|_s\ ;\ k_t=\int_0^t
  \left(\1_{\{X_s=\beta\}}-\1_{\{X_s=\alpha\}}\right)d|k|_s,\label{eq:sys.kps}
\end{align}
where $B$ is a one-dimensional Brownian motion.

If $m$ and $b$ are Lipschitz with respect to the first variable, then
one can simulate this diffusion on $[0,T]$ as follows. Fix $h>0$,
$\bar{\alpha}$ and $\bar{\beta}$ such that
$\alpha<\bar{\alpha}<\bar{\beta}<\beta$, and let $n$ be the
first integer greater than $T/h$. For $\rho\in\{0,1,\ldots,n-1\}$ and
$\rho h< t\leq (\rho+1)h$, let
\begin{align*}
  \tilde{X}_0 & =x,\\
  \tilde{X}_t & =\max[\alpha,\min[\beta,\tilde{X}_{\rho h}+b(\tilde{X}_{\rho
    h},u)(t-\rho h)+\sqrt{2m(\tilde{X}_{\rho h},u)}(B_t-B_{\rho h}) \\
  & \phantom{=\max[\alpha,\min[\beta,\tilde{X}_{\rho h}\ }
  +\mathbf{1}_{\{\tilde{X}_{\rho
      h}<\bar{\alpha}\}}\max(0,A^\rho_t-(\tilde{X}_{\rho h}-\alpha)) \\
  & \phantom{=\max[\alpha,\min[\beta,\tilde{X}_{\rho h}\ }
  -\mathbf{1}_{\{\tilde{X}_{\rho
      h}>\bar{\beta}\}}\max(0,B^\rho_t+(\tilde{X}_{\rho h}-\beta))]],
\end{align*}
where
\begin{align*}
  A^\rho_t & =\sup_{\rho h\leq s\leq t}\left\{-b(\tilde{X}_{\rho h},u)(s-\rho
  h)-\sqrt{2m(\tilde{X}_{\rho h},u)}(B_s-B_{\rho h})\right\}, \\
  B^\rho_t & =\sup_{\rho h\leq s\leq t}\left\{b(\tilde{X}_{\rho h},u)(s-\rho
  h)+\sqrt{2m(\tilde{X}_{\rho h},u)}(B_s-B_{\rho h})\right\}.
\end{align*}
Then, there exists a constant $C$ independent of $h$ such that for
any function $f$ on $[\alpha,\beta]$ with finite variation bounded
by $1$, $\sup_{0\leq t\leq T}|E(f(X_t)-f(\tilde{X}_t))|\leq
C\sqrt{h}$.

In each step of this scheme, one has to simulate simultaneously $B_t$
and $S_t:=\sup_{s\leq t}(aB_s+bs)$ for fixed constants $a,b$. This can be
done as follows (Shepp~\cite{Shepp79}). Let $U$ be a Gaussian centered
random variable with variance $t$, and let $V$ be an exponential
random variable with parameter $1/2t$ independent of $U$. Put
\begin{equation*}
  Y=\frac{1}{2}\left(aU+bt+\sqrt{a^2V+(aU+bt)^2}\right).
\end{equation*}
Then the vectors $(B_t,S_t)$ and $(U,Y)$ have the same distribution.

Note that this scheme can be easily generalized to state spaces of
the form
$(\alpha_1,\beta_1)\times\ldots\times(\alpha_d,\beta_d)\subset\mathbb{R}^d$,
as explained in~\cite{Le95}.

Next, we want to study a particular case in which we can considerably
reduce the complexity of the algorithm. In Section~\ref{sec:mart}, one
needs to compute $I^\delta W*\nu(x,u)$ at some point $(x,u)\in{\cal
  X}\times{\cal U}$ at each time step, which involves a sum over all
individuals in the population. In the case of logistic competition
(linar death rate) where
\begin{equation}
  \label{eq:logist}
  \mu(x,u,r)=\mu_0(x,u)+\mu_1(x,u)r,
\end{equation}
one can use the following algorithm.

Fix a constant $C_{\delta}$ in a similar way as in~(\ref{cdelta1})
such that $\mu_0(x,u)+\lambda^*+\|M^*\|_1\leq C_\delta$ and
$\mu_1(x,u)\|I^{\delta}W\|_{\infty}\leq C_{\delta}$.
Take the Brownian motions $(B^{j,k})_{j,k\in\mathbb{N}}$ and the
random variables $(\theta_k)_{k\in\mathbb{N}}$,
$(V_k)_{k\in\mathbb{N}}$ and $(\tau_k)_{k\in\mathbb{N}}$ as in
Section~\ref{sec:mart}. Set $T_0=0$ and $N_0=N$ (the initial
number of individuals). Assume that
$(T_{k-1},N_{k-1},X_{T_{k-1}},U_{T_{k-1}})$ are given. $N_{k-1}$
is the number of individuals at time $T_{k-1}$. At this time,
their positions and traits are the coordinates of the vectors
$X_{T_{k-1}}=(X^i_{T_{k-1}})_{1\leq i\leq N_{k-1}}$ and
$U_{T_{k-1}}=(U^i_{T_{k-1}})_{1\leq i\leq N_{k-1}}$.  The two
first steps of the algorithm are the same: the new time step is
given by $T_k=T_{k-1}+\tau_k/N_{k-1}(N_{k-1}+1)$ and the motion of
each particle is governed by the SDE with normal
reflection~(\ref{sys.kps}).

The third step deals with the different events that may happen at
time $T_k$. Choose at random one individual $I_k=i$ uniformly
among the $N_{k-1}$ individuals living during the time interval
$[T_{k-1},T_k)$. Its position and trait are
$(X^i_{T_k},U^i_{T_{k-1}})$. The event occurring at time $T_k$ is
decided by comparing $\theta_k$ with constants related to the rate
of each kind of event. The only difference with the algorithm of
Section~\ref{sec:mart} is in the first sub-step, that has to be
divided in two steps as follows:
\begin{itemize}
\item If $0\leq
  \theta_k<\frac{N_{k-1}}{N_{k-1}+1}=:\theta_0^i(X_{T_{k}},U_{T_{k-1}})$,
  then let $j\in\{1,\ldots,N_{k-1}\}$ be such that
  $\frac{j-1}{N_{k-1}+1}\leq \theta_k<\frac{j}{N_{k-1}+1}$. If
  $\theta_k-\frac{j-1}{N_{k-1}+1}<\frac{
    \mu_1(X^{i}_{T_{k}},U^{i}_{T_{k-1}})I^\delta(X^i_{T_k}-X^j_{T_k})
    W(U^{i}_{T_{k-1}}-U^{j}_{T_{k-1}})}{C_\delta}$, then the
  individual $\: i\: $ dies from competition with individual $\:j\:$
  and $N_k=N_{k-1}-1$. Otherwise, nothing happens and $N_k=N_{k-1}$.
\item If $\:\theta_0^i(X_{T_{k}},U_{T_{k-1}})\leq \theta_k\leq
  \theta_0^i(X_{T_{k}},U_{T_{k-1}})+\frac{\mu_0(X^{i}_{T_{k}},U^{i}_{T_{k-1}})}
  {C_\delta(N_{k-1}+1)}=: \theta_1^i(X_{T_{k}},U_{T_{k-1}})$, then the
  individual $\:i\:$ dies by natural death and $N_k=N_{k-1}-1$.
\end{itemize}
The three other sub-steps are the same.

The main difference with the algorithm of Section~\ref{sec:mart}
is that we no longer have to compute
$\sum_{j=1}^{N_{k-1}}I^\delta(X^i_{T_k}-X^j_{T_k})
W(U^{i}_{T_{k-1}}-U^{j}_{T_{k-1}})$ in the first sub-step, but it
suffices to compute $I^\delta(X^i_{T_k}-X^j_{T_k})
W(U^{i}_{T_{k-1}}-U^{j}_{T_{k-1}})$ for chosen $i$ and $j$.
Moreover, we do not need to compute the position of each
individual in the population at each time step. The third step
above only needs to compute the position of at most two particles at
time $T_k$ (the particles numbered $i$ and $j$).

We assume in the following examples a logistic competition of the
form~(\ref{eq:logist}) and a physical space of the form ${\cal
  X}=(\alpha,\beta)$. Our simulations are realized with the previous
algorithm.

\subsection{Example 1. Spatial clustering}
\label{sec:ex1}

We consider here a set of parameters similar to the one
of~\cite{DD03} and~\cite{Prevost:04}, in which, for each spatial
position $x$, the growth rate is maximal for the trait value
$u=x$. This can represent the effect of a gradual spatial
distribution of different resources, involving a gradual
distribution of traits. For example, for some bird species, a
linearly spatially varying seed size determines a linear variation
of the beak sizes (Grant and Grant \cite{GG02}).
\begin{gather*}
  {\cal X}=(0,1),\quad{\cal U}=[0,1],\quad m(x,u)\equiv m,\quad
  b(x,u)\equiv 0,\\ \lambda(x,u)=2-20(x-u)^2 \mbox{\ if\ }|x-u|\leq
  1/\sqrt{10};\ 0 \mbox{\ otherwise,}\\ \mu_N(x,u,r)=1+\frac{r}{N},\quad
  I^\delta(y)=C_\delta\1_{\{|y|\leq\delta\}},\quad W(v)\equiv 1.
\end{gather*}
Moreover, $M(x,u,v)= 0.1 \times k_s(u,v)$ where $0.1$ is the
mutation rate and $k_s(u,v)$ is the probability density of a
Gaussian random variable with mean $u$ and standard deviation $s$
conditioned on staying in ${\cal U}=[0,1]$. Therefore, we have
four free parameters in this model, $m, \delta, s$ and the
population size $N$. The initial population in our simulations is
composed of $N$ individuals at $(0.5,0.5)$.

\begin{figure}
  \centering
  \mbox{\subfigure[$N=3000, s=0.01, m=0.01, \delta=0.3$.]%
{\epsfig{bbllx=50pt,bblly=50pt,bburx=554pt,bbury=770pt,%
figure=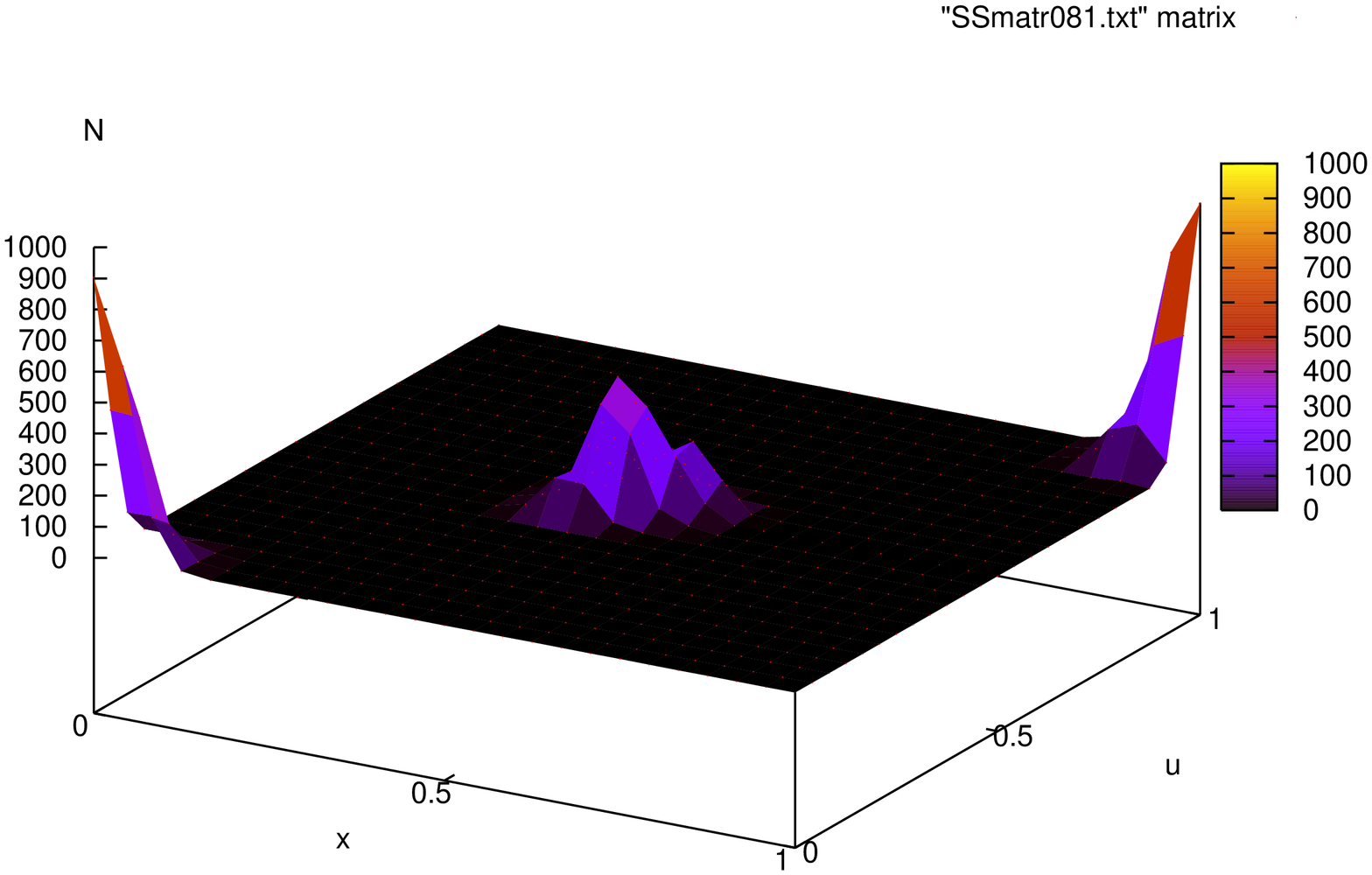, width=.38\textwidth,angle=270}}\quad
    \subfigure[$N=3000, s=0.01, m=0.01, \delta=0.1$.]%
{\epsfig{bbllx=50pt,bblly=50pt,bburx=554pt,bbury=770pt,%
figure=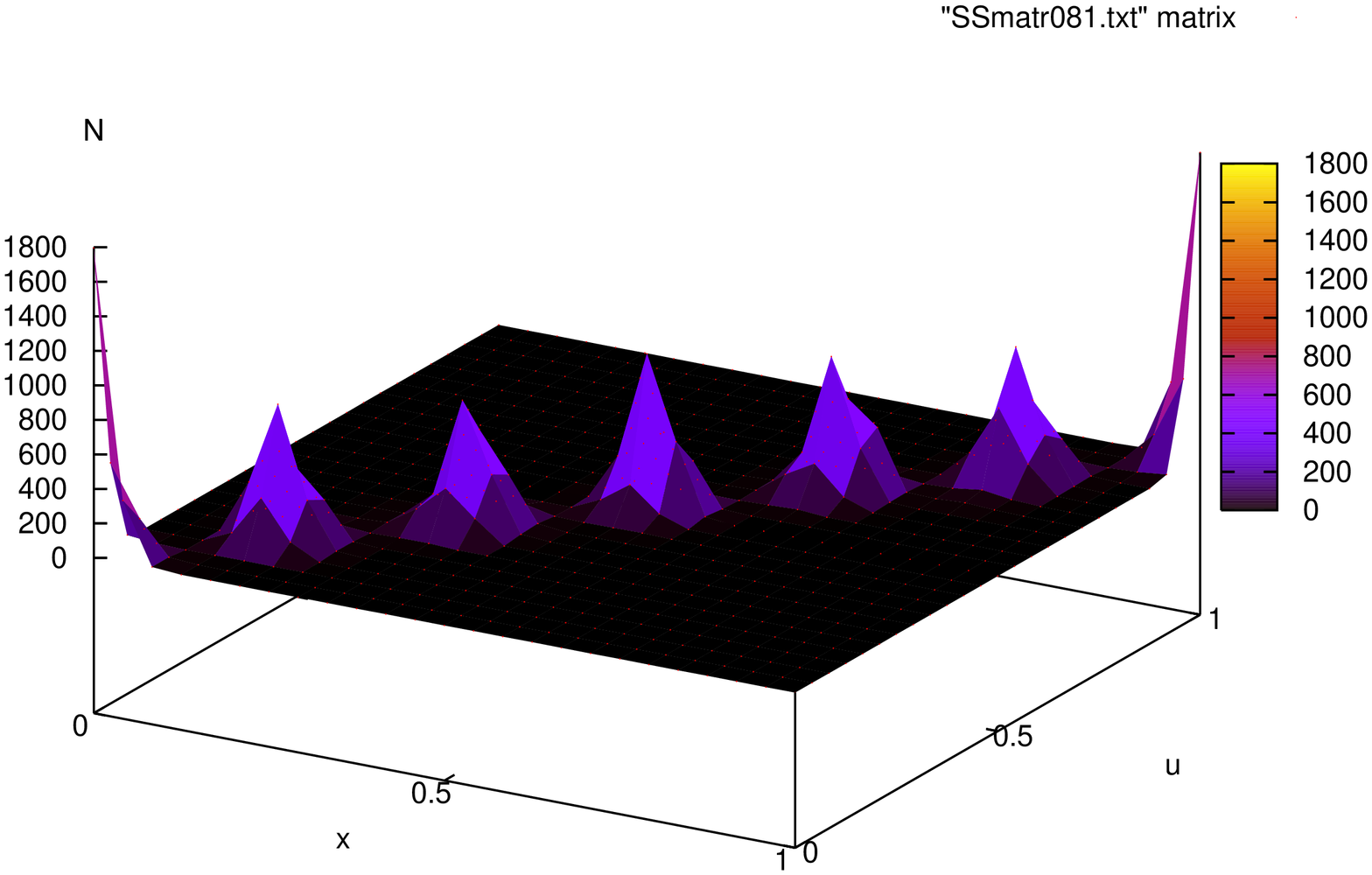, width=.38\textwidth,angle=270}}} \\
  \mbox{\subfigure[$N=3000, s=0.01, m=0.01, \delta=0.03$.]%
{\epsfig{bbllx=50pt,bblly=50pt,bburx=554pt,bbury=770pt,%
figure=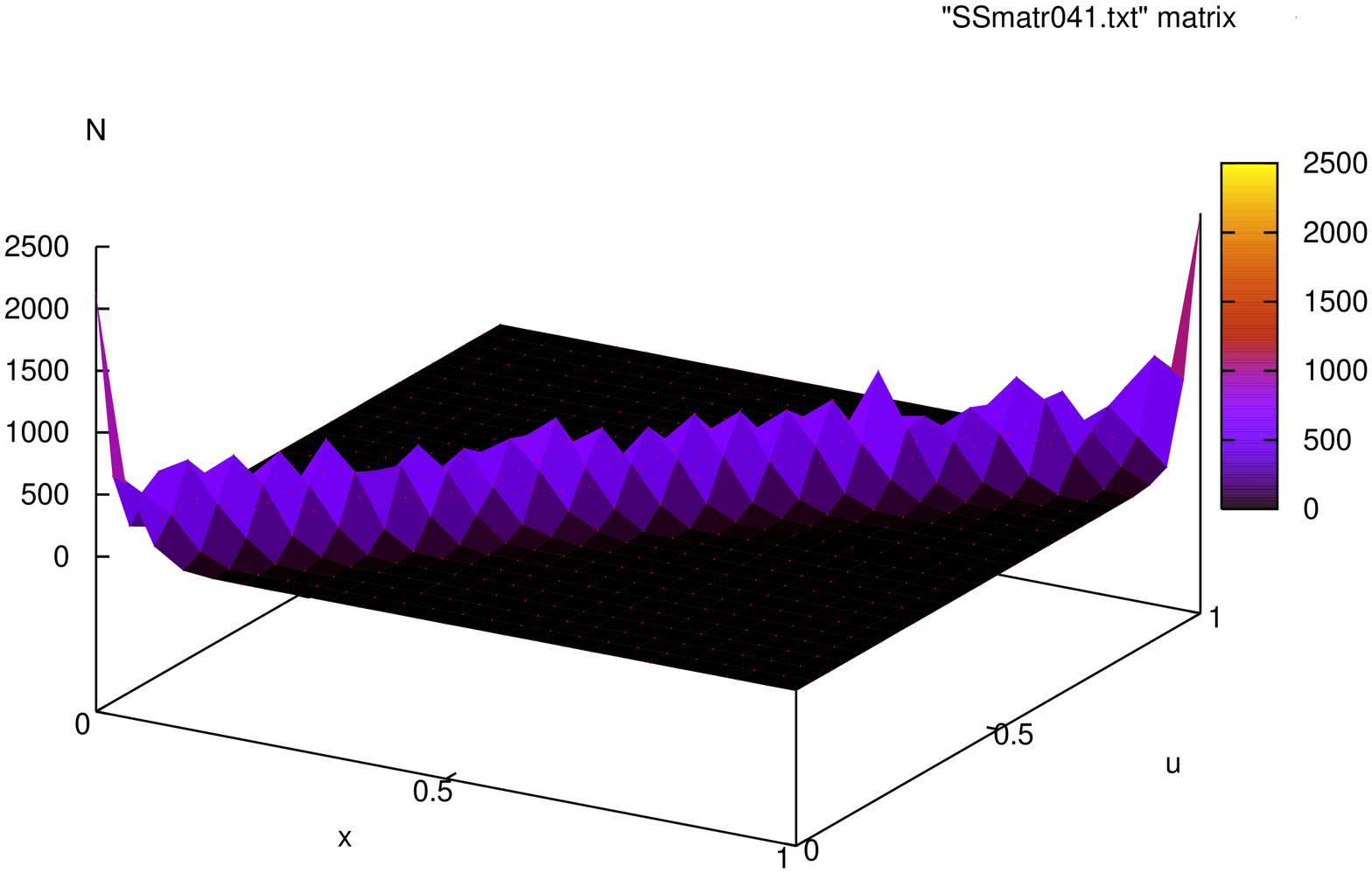, width=.38\textwidth,angle=270}}\quad
    \subfigure[$N=50, s=0.01, m=0.01, \delta=0.1$.]%
{\epsfig{bbllx=50pt,bblly=50pt,bburx=554pt,bbury=770pt,%
figure=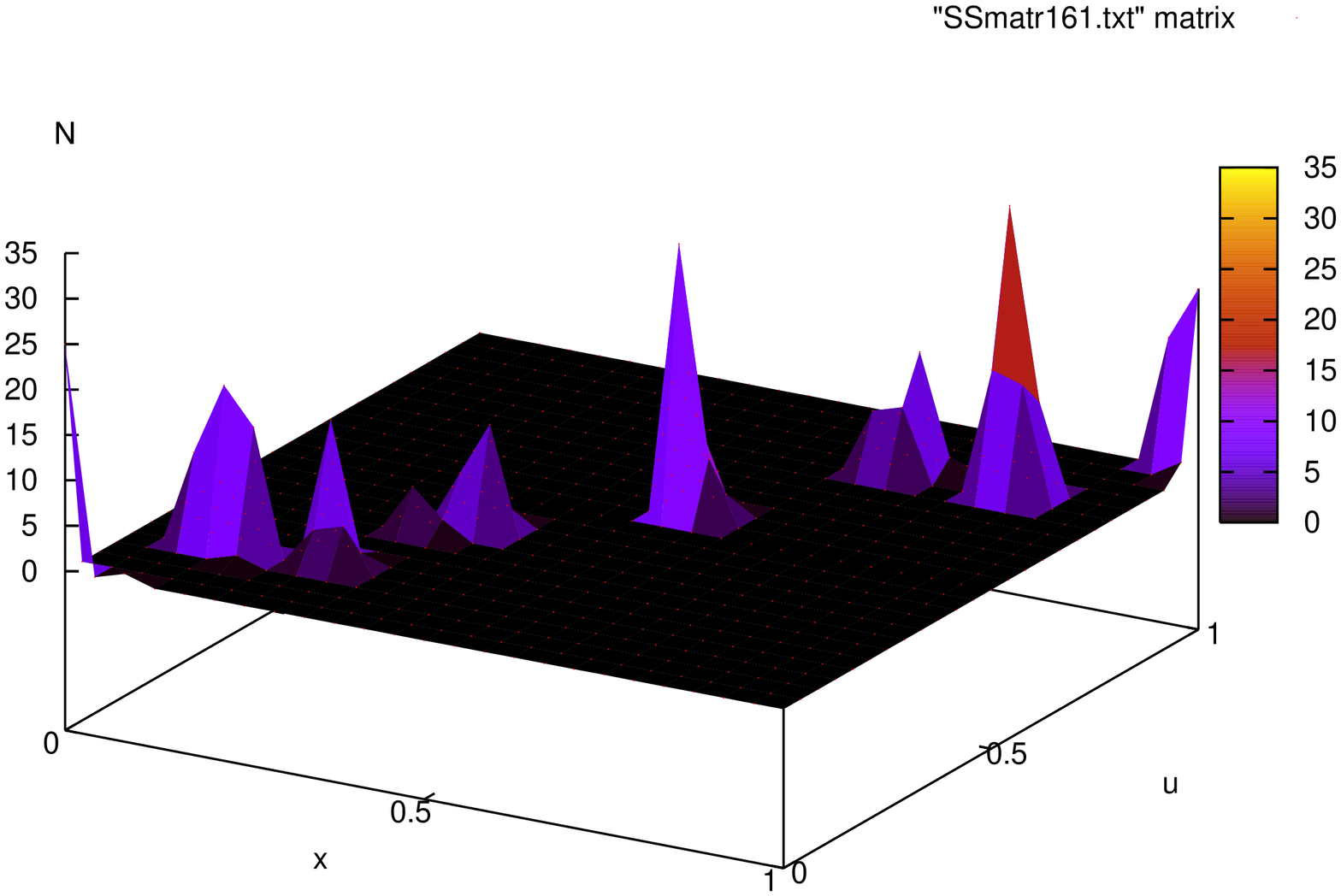, width=.38\textwidth,angle=270}}} \\
  \mbox{\subfigure[$N=3000, s=0.003, m=0.03, \delta=0.1$.]%
{\epsfig{bbllx=50pt,bblly=50pt,bburx=554pt,bbury=770pt,%
figure=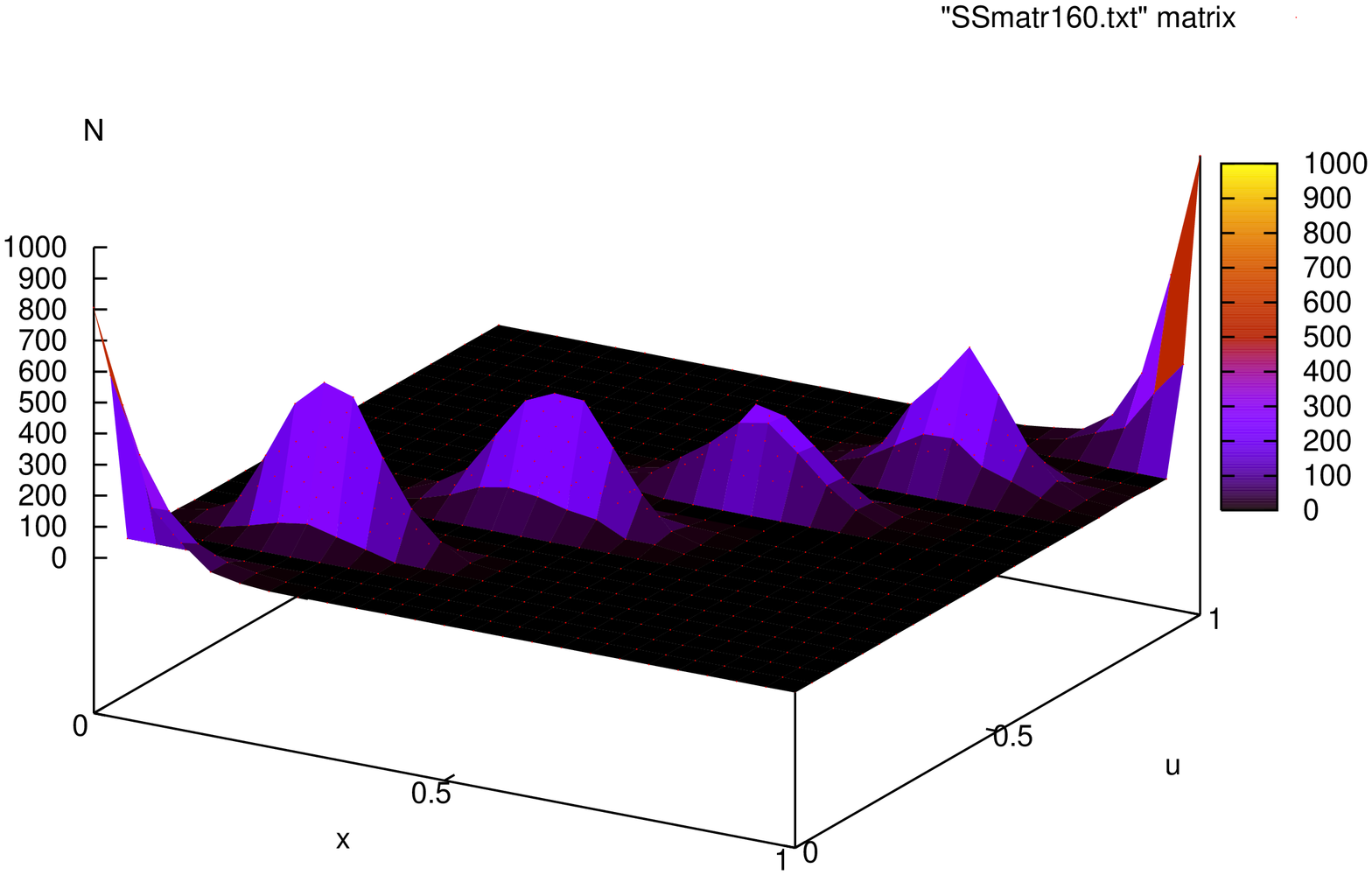, width=.38\textwidth,angle=270}}\quad
    \subfigure[$N=3000, s=0.03, m=0.03, \delta=0.1$.]%
{\epsfig{bbllx=50pt,bblly=50pt,bburx=554pt,bbury=770pt,%
figure=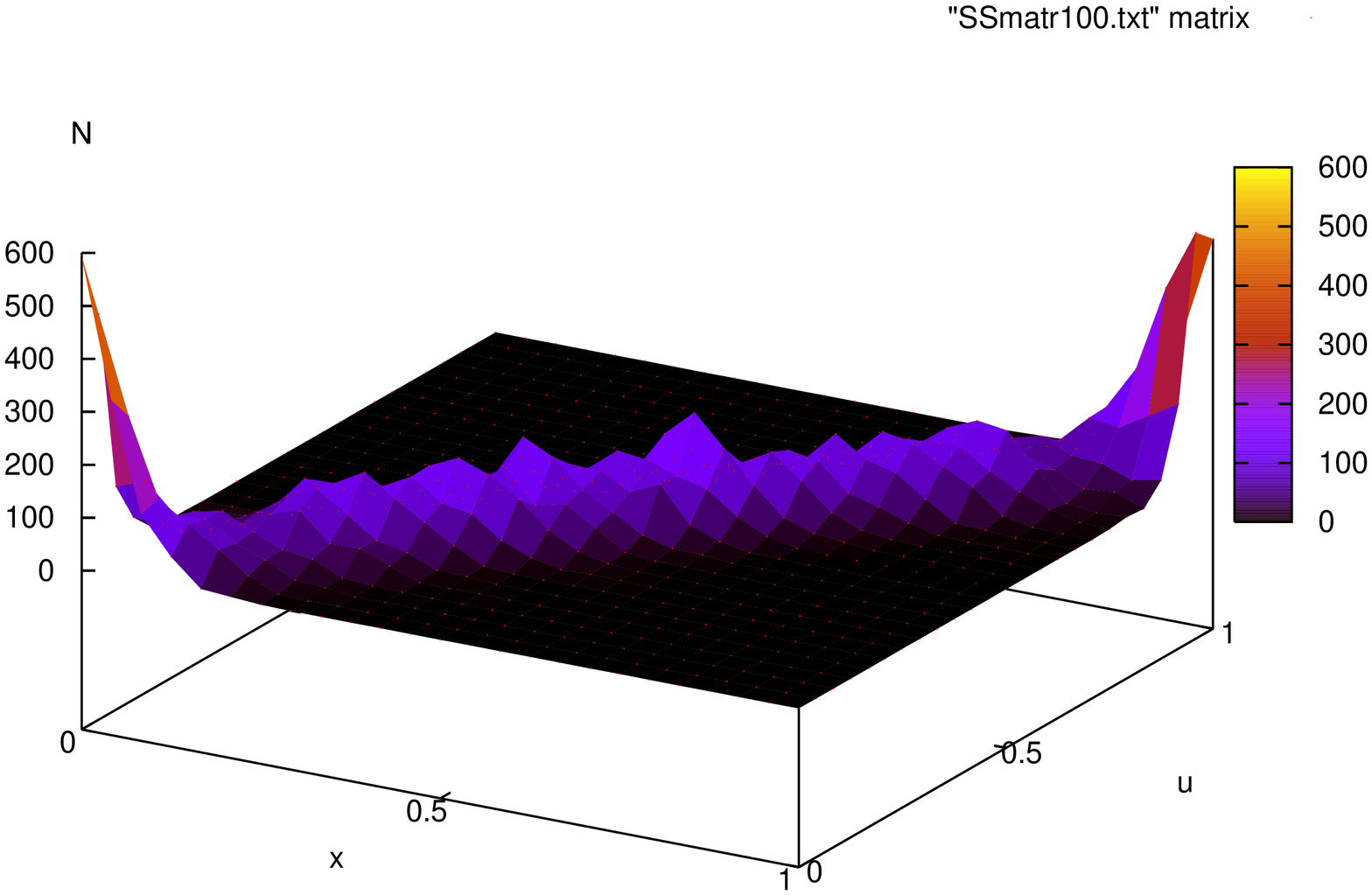, width=.38\textwidth,angle=270}}}
  \caption{Simulations of example 1 for various parameters. All of
    them are taken at time 4000.}
  \label{fig:ex1}
\end{figure}

The simulations of this model show, as in \cite{Prevost:04}, that the
invasion of space occurs along the diagonal $x=u$, and, as
in~\cite{DD03}, that speciation (stable coexistence of several
sub-populations with different typical traits) may occur in this
model, accompanied with a spatial specialization. Several different
population clusters may coexist at different position, with trait
values located around the corresponding optimal traits.  We have
investigated in our simulations the effect of the different parameters
on the clustering and polymorphism of the population. We give pictures
of the seemingly stable state of the population (Fig.~\ref{fig:ex1}).
Our first general observation is that the clusters are more
concentrated at the boundary of the domain. Indeed, the reflected
diffusion governing the motion of individuals is not isotropic close
to this boundary, so that the population density is bigger.
% The qualitative results of this example are also valid for the other
% examples of this section, as well as for many other models we tried.

On the one hand, we have investigated the effect of the
interaction range $\delta$. The main result is that the existence,
the number of clusters and the distance between clusters strongly
depend on the interaction range $\delta$. As shown in
Fig.~\ref{fig:ex1}(a--c), the number of clusters decreases with
$\delta$ and the distance between two population peaks is roughly
$2\delta$, which is exactly the width of the interaction interval.
The emergence of population clusters is mainly a consequence of
local births (\cite{YRS01}). Indeed, since the progeny of an
individual is close to its original location, each individual's progeny can
create a colony with stable position on short timescales. Once an
individual is at a distance greater than $\delta$ from the main
part of the population, it experiences very little competition and
it can create a new colony. When several colonies appear, they
organize in a way to minimize the competition between them and to
maximize the growth rate.
%The approximate distance of $2\delta$
%between two clusters implies that no individual inside a cluster
%of width $\delta$ feels the competition from individuals in the
%center of the other clusters.

If $\delta$ is sufficiently small, we observe a flat distribution
of the population (Fig.~\ref{fig:ex1}(c)), and thus a qualitative
difference with respect to cases (a--b). As an explanation,
decreasing $\delta$ increases the number of clusters and the width
of a cluster increases with the speed of dispersal $m$ and the
range of mutation $s$.
 Then no distinct colony can be observed for sufficiently small
$\delta$, and fixed $m$ and $s$.

We also investigated the effect of the population size $N$. It
appears that this parameter has very little qualitative effect on
the clustering of the population, except for small $N$
(Fig.~\ref{fig:ex1}(d)), where the width of each clusters is
reduced, and we observe much more fluctuations in the population
distribution. However, we still can observe a similar pattern of
population clusters than in Fig.~\ref{fig:ex1}(b).

On the other hand, we also studied the effect of the diffusion
coefficient $m$ and the mutation range $s$. Comparing
Fig.~\ref{fig:ex1}(b) and (f), we observe that too large $s$ and
$m$ induce the same flat distribution as for small $\delta$. This
confirms that the clusters pattern depends mainly on the balance
between $m$ and $s$, and $\delta$. In Fig.~\ref{fig:ex1}(f),
quick movements mix the population so that no spatial structure
can appear.

Finally, we also studied the relative effect of $s$ and $m$ in the
appearance of spatial or phenotypic structure. As shown in
Fig.~\ref{fig:ex1}(e), small $s$ can induce a differentiation over the
phenotype space ${\cal U}$ even when $m$ is large enough to have a
flat distribution over space ${\cal X}$ (compare with
Fig.~\ref{fig:ex1}(f)). Fig.~\ref{fig:ex1}(e) can be seen as an
intermediate state between Fig.~\ref{fig:ex1}(a) and
Fig.~\ref{fig:ex1}(f). When $m$ is reduced instead of $s$, a reversed
pattern can be observed.

\subsection{Example 2. The role of spatial competition for clustering}
\label{sec:ex2}

As we have seen above, the balance between the spatial competition
range $\delta$ and the diffusion parameters $s$ and $m$ has an
important effect on the clustering of the population. Here we want
to address  the balance between the range of competition and the
growth rate. For this purpose, we consider the following model,
inspired by the adaptive dynamics model of~\cite{DD99}:
\begin{gather*}
  {\cal X}=(-1,1),\quad{\cal U}=[0,2],\quad m(x,u)\equiv m,\quad
  b(x,u)\equiv 0,\\ \lambda(x,u)=\exp(-x^2/2\rho^2),\quad
  \mu_N(x,u,r)=1+\frac{r}{N},\\
  I^\delta(y)=C_\delta\exp(-y^2/2\delta^2),\quad W(v)=\exp(-v^2/0.02).
\end{gather*}
and the same mutation kernel as above. This example has five free
parameters $m, \delta, s, N$,  and $\rho$, which represents the
width of the space region with significant growth rate (namely, a
parameter describing the width of the space region with high
concentration of resources). The initial population in our
simulations is composed of $N$ individuals at $(0,1)$. Observe
that in this example, the trait has no effect on the growth rate,
so that the trait structure is neutral (all individual's
parameters are equal, independent of the trait, in absence of
interaction).

\begin{figure}
  \centering
  \mbox{\subfigure[$\delta=0.9$, $\rho=1$.]%
{\epsfig{bbllx=50pt,bblly=50pt,bburx=554pt,bbury=770pt,%
figure=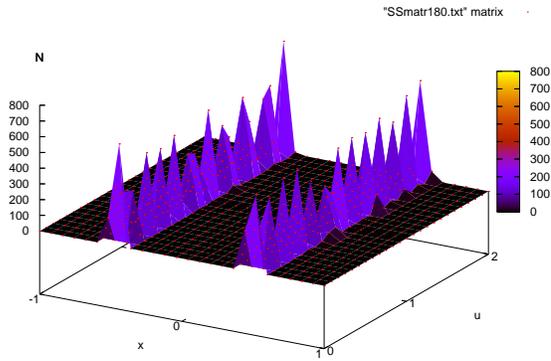, width=.38\textwidth,angle=270}}\quad
    \subfigure[$\delta=1.1$, $\rho=1$.]%
{\epsfig{bbllx=50pt,bblly=50pt,bburx=554pt,bbury=770pt,%
figure=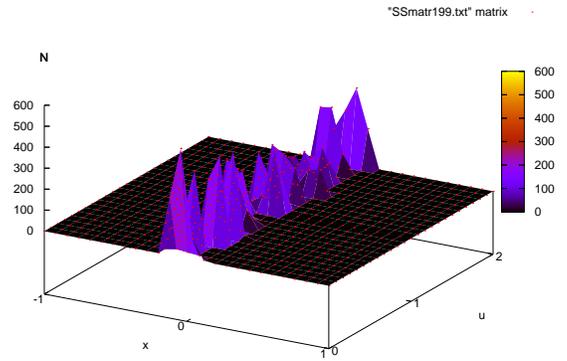, width=.38\textwidth,angle=270}}} \\
  \mbox{\subfigure[$\delta=0.5$, $\rho=1$.]%
{\epsfig{bbllx=50pt,bblly=50pt,bburx=554pt,bbury=770pt,%
figure=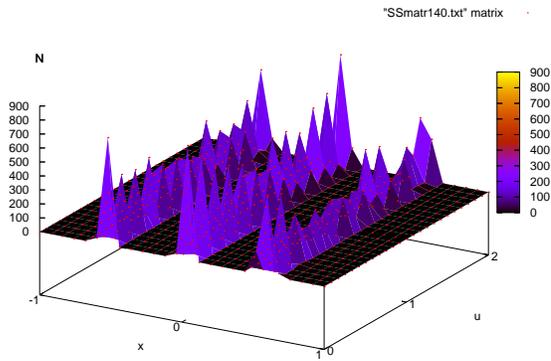, width=.38\textwidth,angle=270}}\quad
    \subfigure[$\delta=0.1$, $\rho=1$.]%
{\epsfig{bbllx=50pt,bblly=50pt,bburx=554pt,bbury=770pt,%
figure=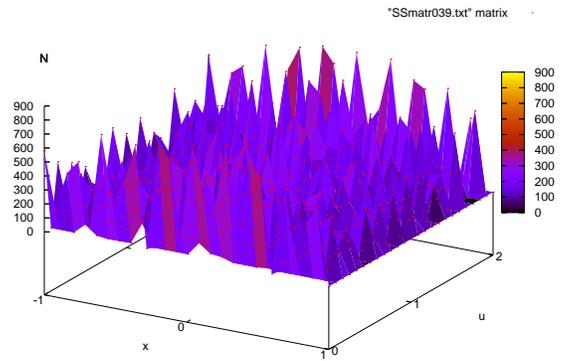, width=.38\textwidth,angle=270}}}
\caption{Simulations of example 2. Neutral case. All the figures are
  taken at time 5000, except the last one, taken at time 10000 (more
  time is needed to fill the whole space). In all the simulations,
  $N=1000$, $s=0.003$ and $m=0.003$.}
  \label{fig:ex2-init}
\end{figure}

Remark that if we consider the space ${\cal X}$ as a trait space,
this model is similar to the one of~\cite{DD99}. In particular,
the biological theory of adaptive dynamics (\cite{Gal97}) suggests
that evolutionary branching, i.e.\ the split of the population
into two sub-populations with different traits stably coexisting,
translating in our model into spatial clustering, occurs when the
range of interaction ($\delta$ in our case) is smaller than the
range of the growth rate ($\rho$ in our case). This is illustrated
by Fig.~\ref{fig:ex2-init}(a) and (b), where, in (a),
$\delta<\rho$ and the population stabilizes around two distinct
positions (branching occurs) and in (b),  $\delta>\rho$ and the
population stabilizes around position 0 (there is no branching).

Figures (c) and (d) prove that other phase transitions occur for
smaller $\delta$, leading to the coexistence of three clusters or
more. As in example 1, we notice in Figure (d) that very small
$\delta$ leads to a distribution without distinct clusters.
\\

It is possible to add some phenotypic structure to this example by
assuming that the growth rate depends on the trait $u$, in a way such
that spatial branching occurs for some traits but not for others,
according to the above branching criterion. We take the same
parameters, except for the birth rate, which has the following form.
\begin{gather*}
  \lambda(x,u)=\exp(-x^2/2(u+0.1)).
\end{gather*}
The parameter $\rho$ is then replaced by $\sqrt{u+0.1}$, so that
branching occurs if $\sqrt{u+0.1}>\delta$.

\begin{figure}
  \centering
  \mbox{\subfigure[$t=10000$.]%
{\epsfig{bbllx=50pt,bblly=50pt,bburx=554pt,bbury=770pt,%
figure=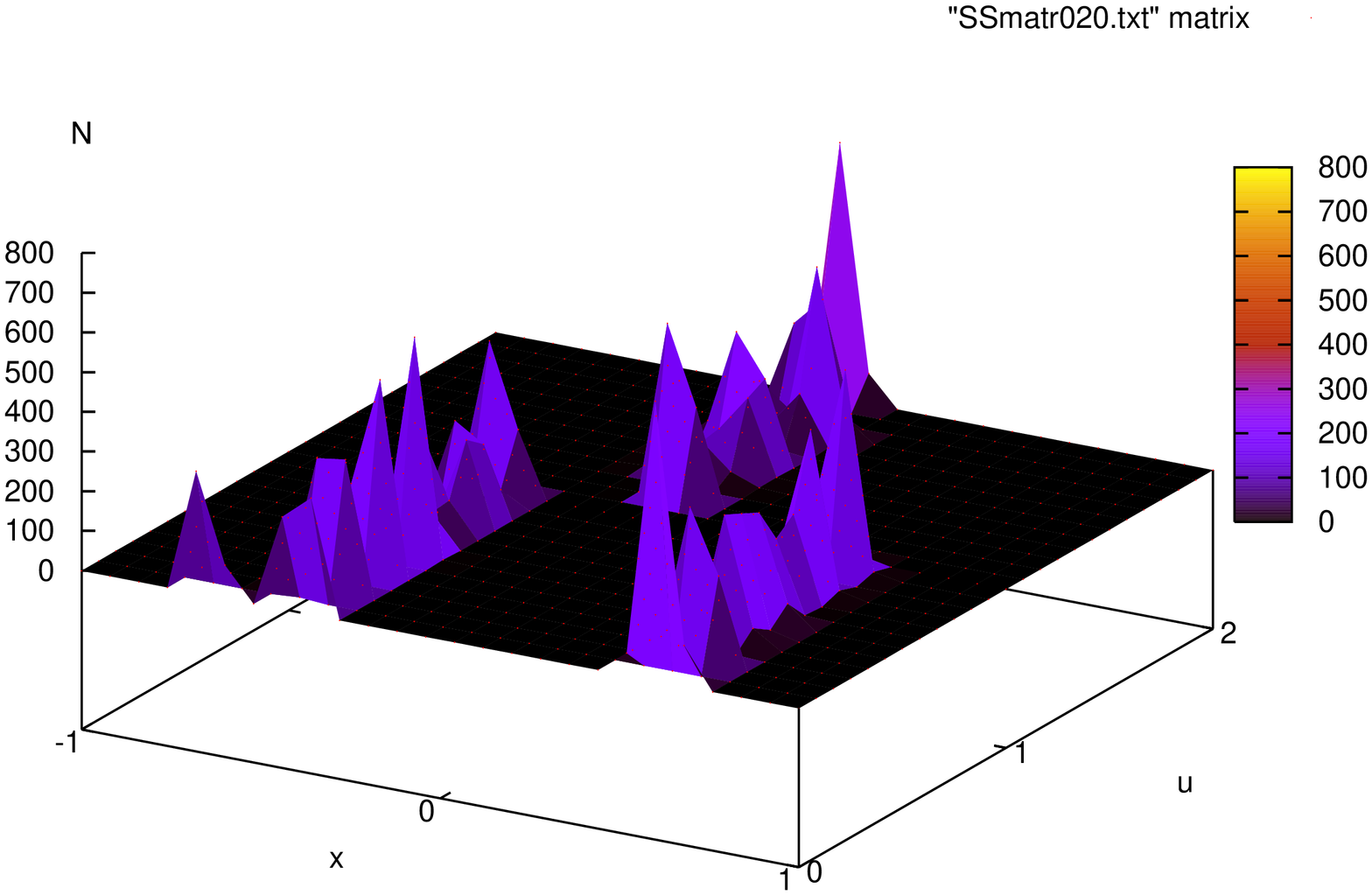, width=.38\textwidth,angle=270}}\quad
    \subfigure[$t=20000$.]%
{\epsfig{bbllx=50pt,bblly=50pt,bburx=554pt,bbury=770pt,%
figure=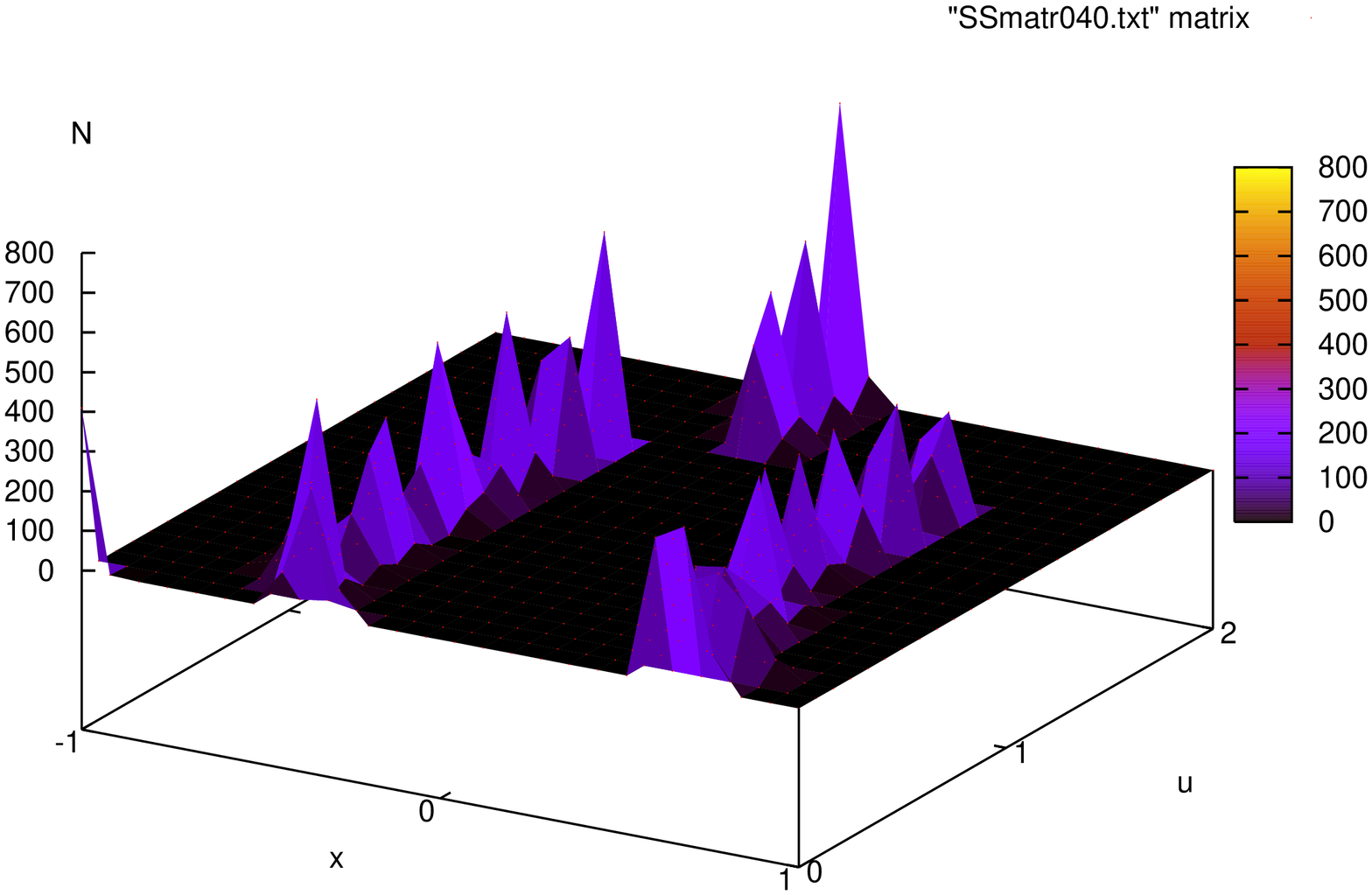, width=.38\textwidth,angle=270}}} \\
  \mbox{\subfigure[$t=45000$.]%
{\epsfig{bbllx=50pt,bblly=50pt,bburx=554pt,bbury=770pt,%
figure=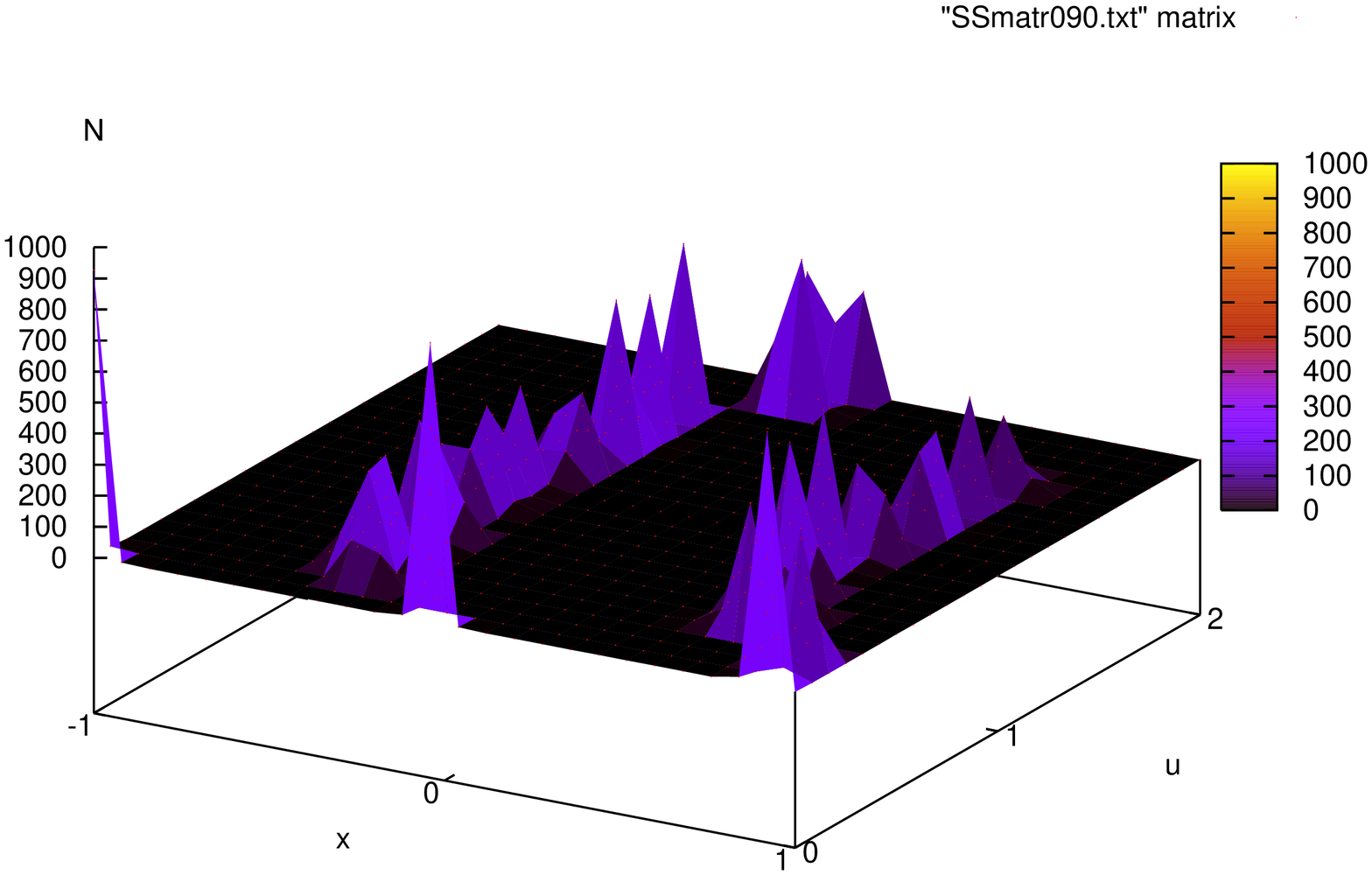, width=.38\textwidth,angle=270}}\quad
    \subfigure[$t=80000$.]%
{\epsfig{bbllx=50pt,bblly=50pt,bburx=554pt,bbury=770pt,%
figure=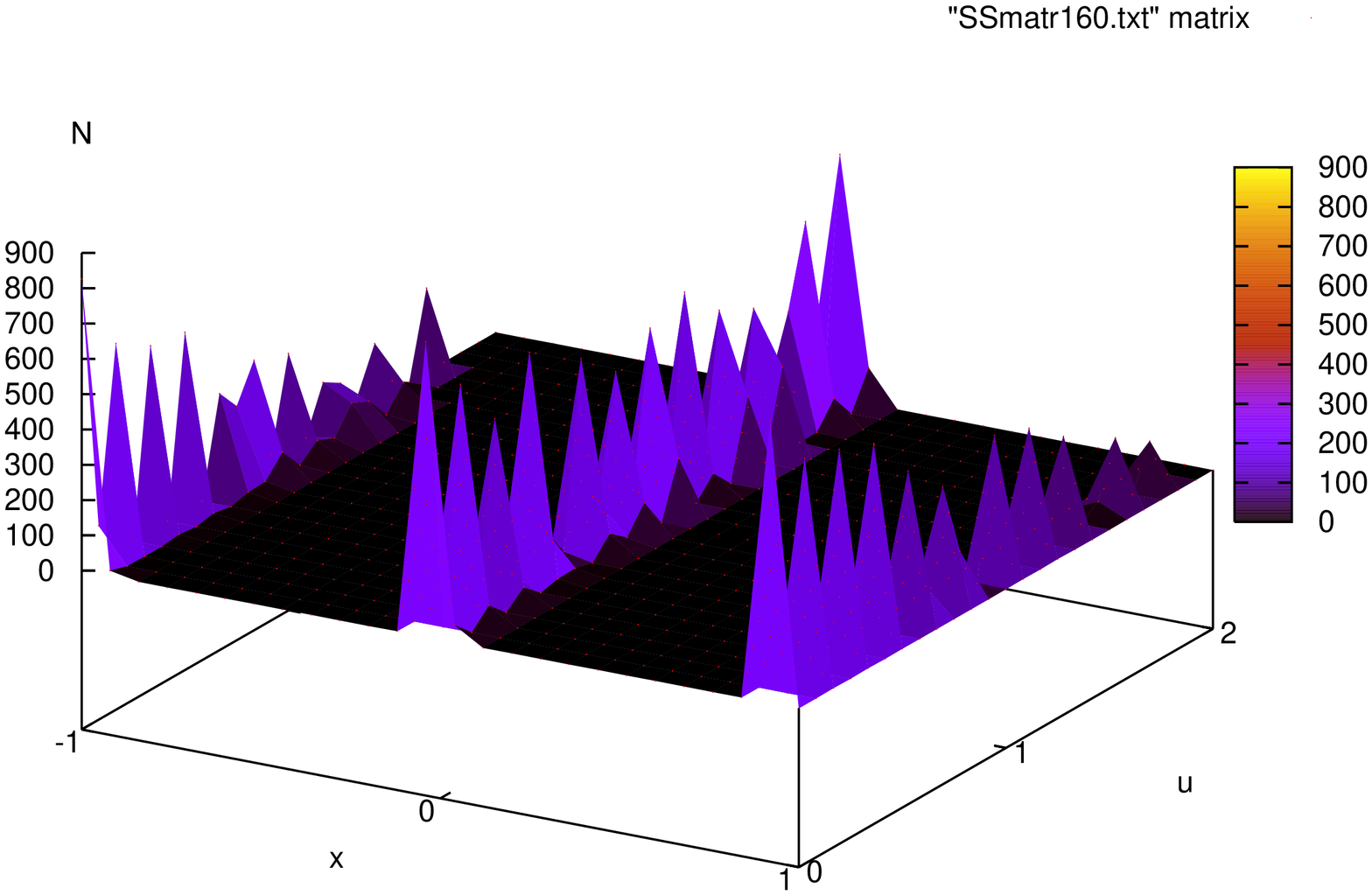, width=.38\textwidth,angle=270}}}
  \caption{Simulations of example 2. Trait-dependent case. In this
    simulation, $N=1000, s=0.003, m=0.003$ and $\delta=1$.}
  \label{fig:ex2}
\end{figure}

This is what happens actually for small times
(Fig.~\ref{fig:ex2}(a)), but after a longer time
(Fig.~\ref{fig:ex2}(b,c)), the two clusters spread over the trait
space because of the mutations. Eventually, if we let time go on,
we actually observe the appearance and the spread of three spatial
clusters, colonizing all the trait space (Fig.~\ref{fig:ex2}(d)).

\subsection{Example 3. Invasion and evolution of migration speed}
\label{sec:invas}

Here, we investigate a model describing the invasion of a species
with evolving dispersal speed (as in \cite{Desville:04}). This can
model phenomena such as the invasion of Australia by cane toads,
for which an adaptation to high invasion speeds has been recently
detected (Phillips et al. \cite{Pal06}). The parameters are as
follows.
\begin{gather*}
  {\cal X}=(-1,1),\quad{\cal U}=[0,3],\quad m(x,u)\equiv m(u+0.1),\quad
  b(x,u)\equiv 0,\\ \lambda(x,u)=1,\quad
  \mu_N(x,u,r)=1+\frac{r}{N},\\
  I^\delta(y)=C_\delta\1_{\{|y|\leq\delta\}},\quad W(v)=\exp(-10v^2).
\end{gather*}
and the same mutation kernel as above.  Here we study invasion
into an homogeneous space ($\lambda$ is constant). The diffusion
rate $m$ is proportional (up to a constant) to the trait $u$.
Thus, individuals with large $u$ move fast. The trait $u$ can be a
morphological trait responsible for the speed of dispersal (e.g.\
the length of legs for toads, \cite{Pal06}). Space competition
occurs between individuals within a distance $\delta$, and the
kernel $W$ models competition between close traits.  This example
has four free parameters, the diffusion coefficient $m$, the
interaction range $\delta$, the standard deviation of mutations
$s$ and the population size $N$.

\begin{figure}
  \centering
  \mbox{\subfigure[$t=25$.]%
{\epsfig{bbllx=50pt,bblly=50pt,bburx=554pt,bbury=770pt,%
figure=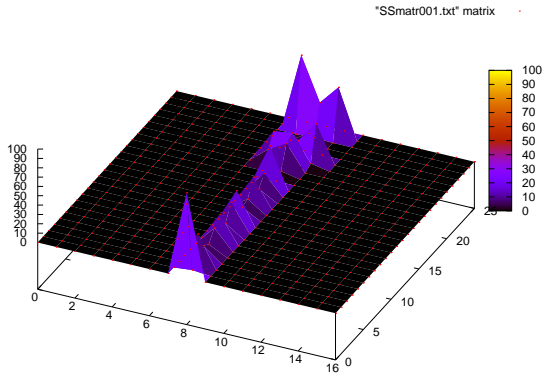, width=.38\textwidth,angle=270}}\quad
    \subfigure[$t=125$.]%
{\epsfig{bbllx=50pt,bblly=50pt,bburx=554pt,bbury=770pt,%
figure=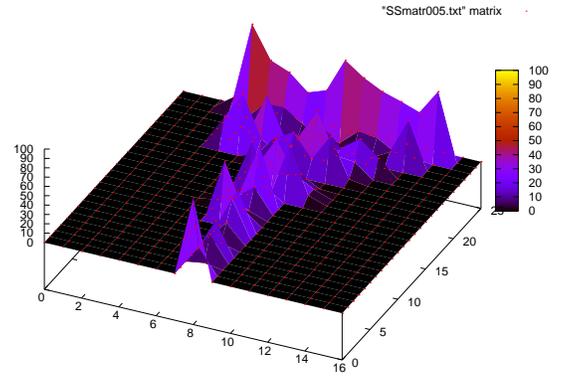, width=.38\textwidth,angle=270}}} \\
  \mbox{\subfigure[$t=250$.]%
{\epsfig{bbllx=50pt,bblly=50pt,bburx=554pt,bbury=770pt,%
figure=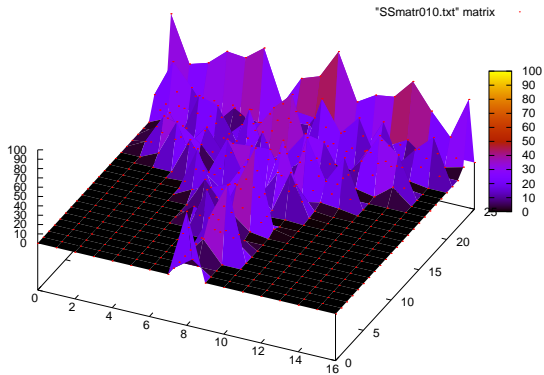, width=.38\textwidth,angle=270}}\quad
    \subfigure[$t=500$.]%
{\epsfig{bbllx=50pt,bblly=50pt,bburx=554pt,bbury=770pt,%
figure=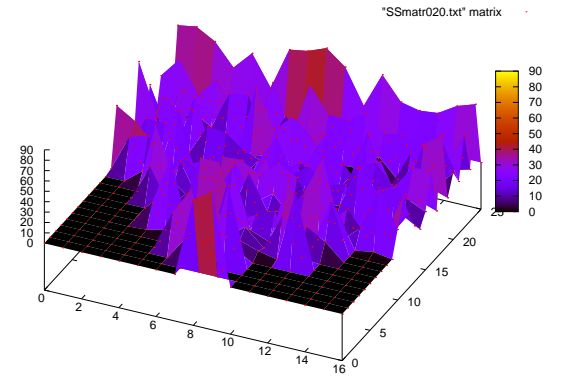, width=.38\textwidth,angle=270}}}
  \caption{Simulations of example 3. The parameters are $N=100, s=0.03,
  m=0.003$ and $\delta=0.1$. The initial condition is composed of $N$
  individuals located at 0 and with trait values $3i/N$ for $1\leq
  i\leq N$.}
  \label{fig:ex8}
\end{figure}

\begin{figure}
  \centering
  \mbox{\subfigure[$t=500$.]%
{\epsfig{bbllx=50pt,bblly=50pt,bburx=554pt,bbury=770pt,%
figure=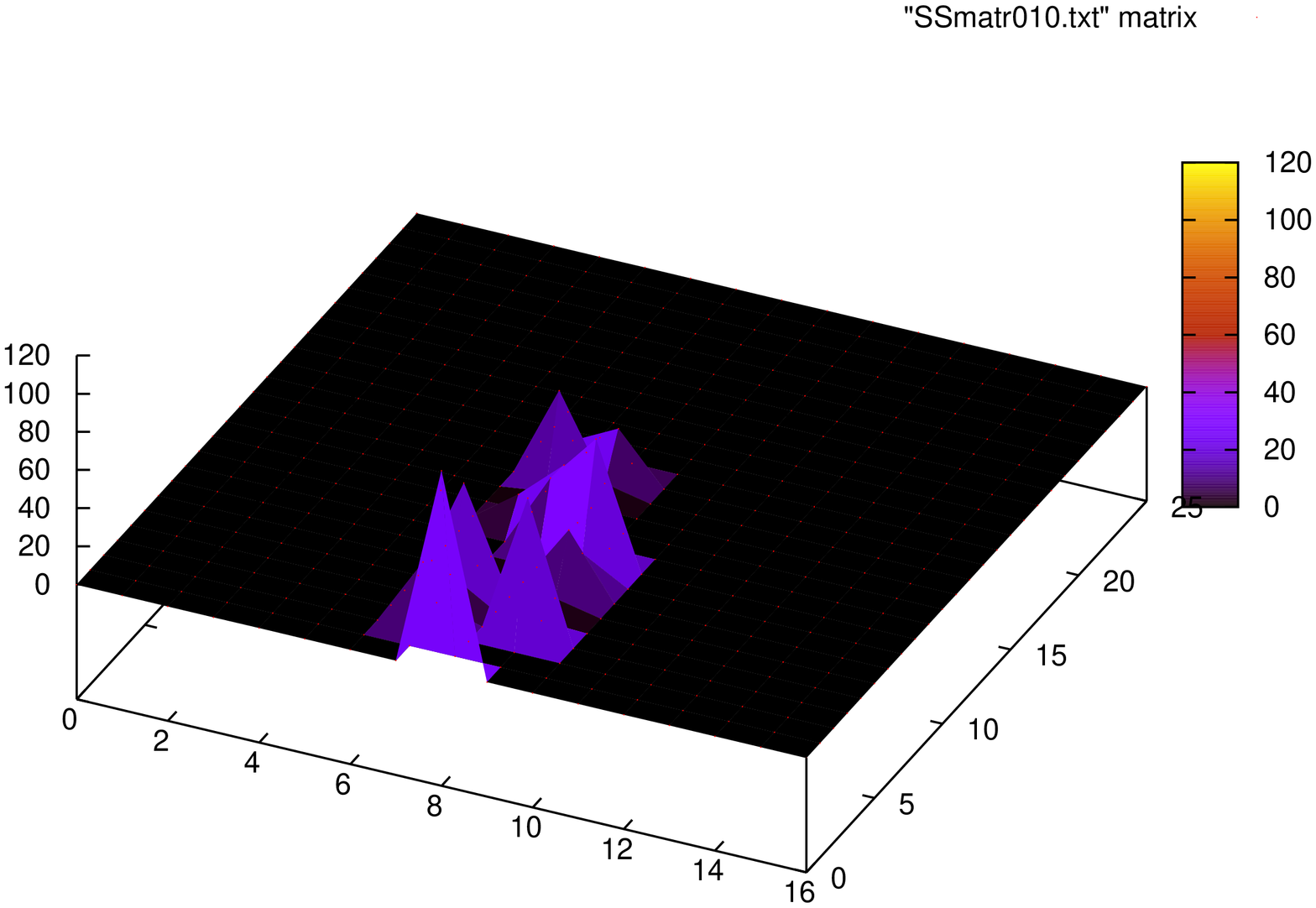, width=.38\textwidth,angle=270}}\quad
    \subfigure[$t=750$.]%
{\epsfig{bbllx=50pt,bblly=50pt,bburx=554pt,bbury=770pt,%
figure=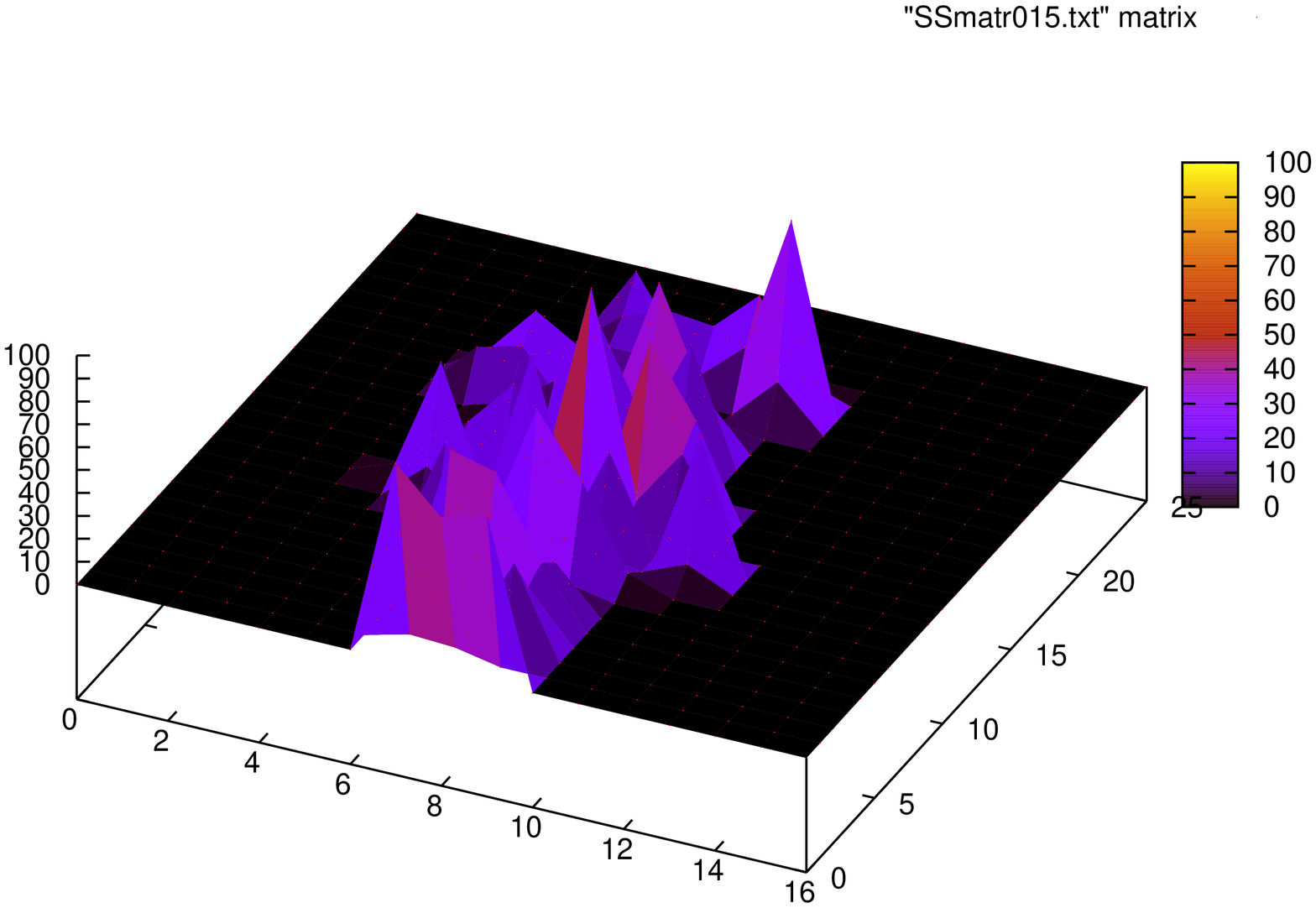, width=.38\textwidth,angle=270}}} \\
  \mbox{\subfigure[$t=850$.]%
{\epsfig{bbllx=50pt,bblly=50pt,bburx=554pt,bbury=770pt,%
figure=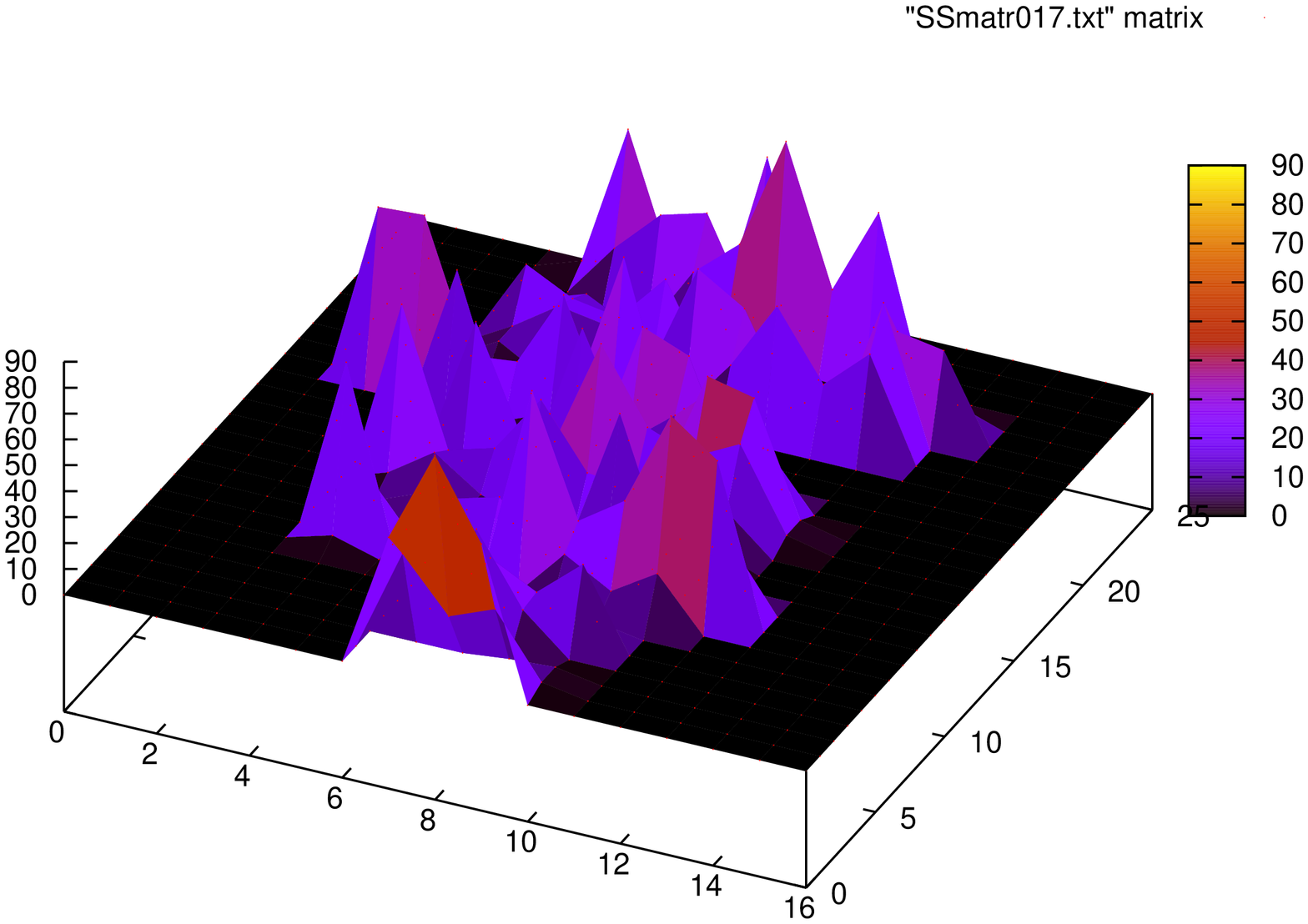, width=.38\textwidth,angle=270}}\quad
    \subfigure[$t=1000$.]%
{\epsfig{bbllx=50pt,bblly=50pt,bburx=554pt,bbury=770pt,%
figure=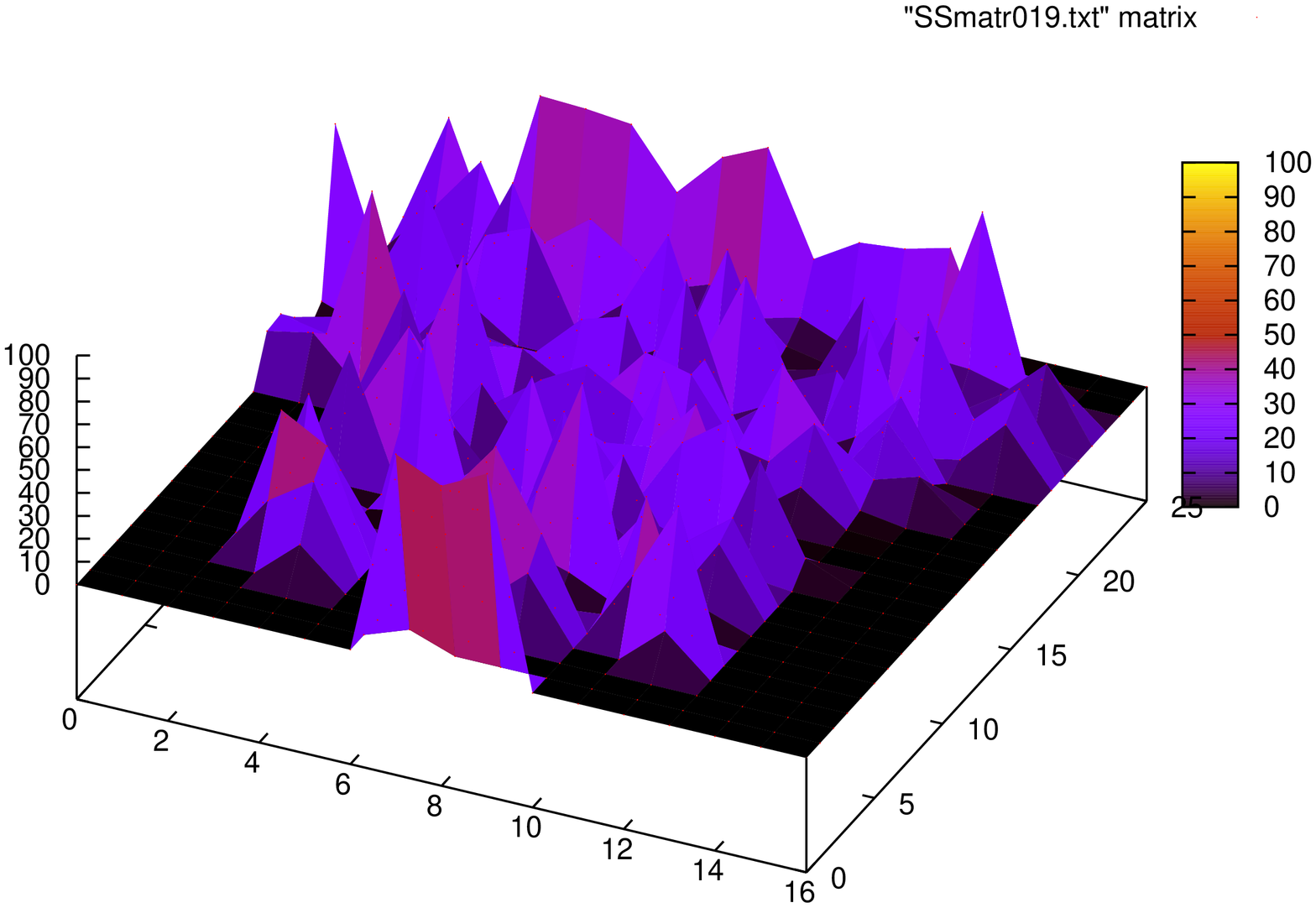, width=.38\textwidth,angle=270}}}
\caption{Simulations of example 3. The parameters are $N=100, s=0.03,
  m=0.003$ and $\delta=0.1$. The initial condition is composed of $N$
  individuals located at $(0,0)$.}
  \label{fig:ex7}
\end{figure}

In Fig.~\ref{fig:ex8} and~\ref{fig:ex7}, we present two extreme cases
with respect to the initial trait distribution, but with identical
parameters. In the first one, all individuals are at (physical)
position 0, and with traits regularly distributed in ${\cal U}=[0,3]$.
In the second one, all individuals are initially located at a single
point $(0,0)$.

In both figures, we observe a triangular invasion pattern
indicating that the invasion front is composed of faster
individuals. In Fig.~\ref{fig:ex7}, we also observe a simultaneous
invasion in physical and trait spaces, and a slower spread of the
population. This can be explained by the fact that the population,
initially composed of slow individuals, has first to colonize the
trait space before invading the physical space.  Because of the
progressive appearance of larger traits, the invasion speed
increases over time (compare the different time values in
Fig.~\ref{fig:ex7}).

\begin{figure}
  \centering
  \mbox{\subfigure[$t=5$.]%
{\epsfig{bbllx=50pt,bblly=50pt,bburx=554pt,bbury=770pt,%
figure=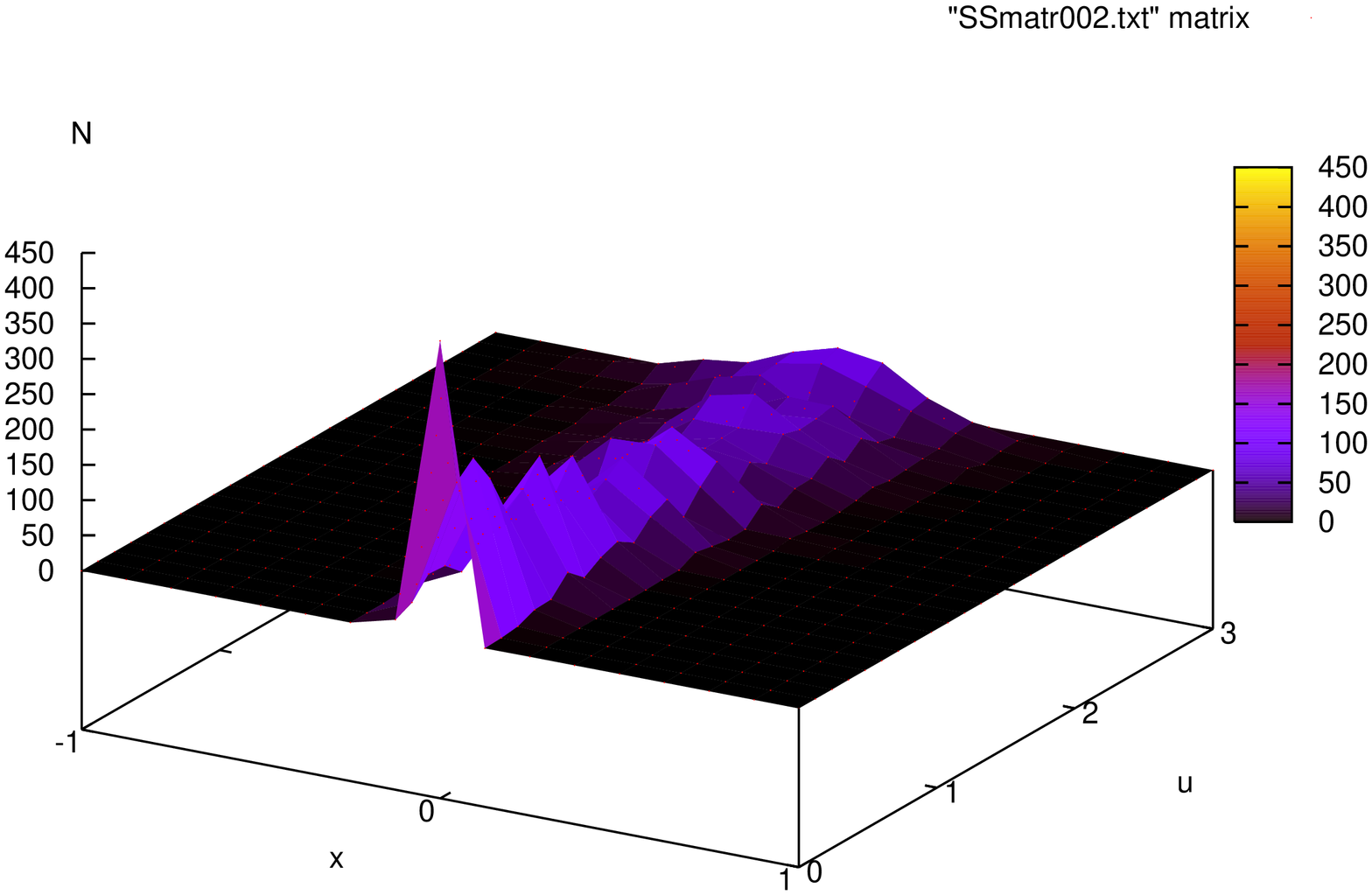, width=.38\textwidth,angle=270}}\quad
    \subfigure[$t=10$.]%
{\epsfig{bbllx=50pt,bblly=50pt,bburx=554pt,bbury=770pt,%
figure=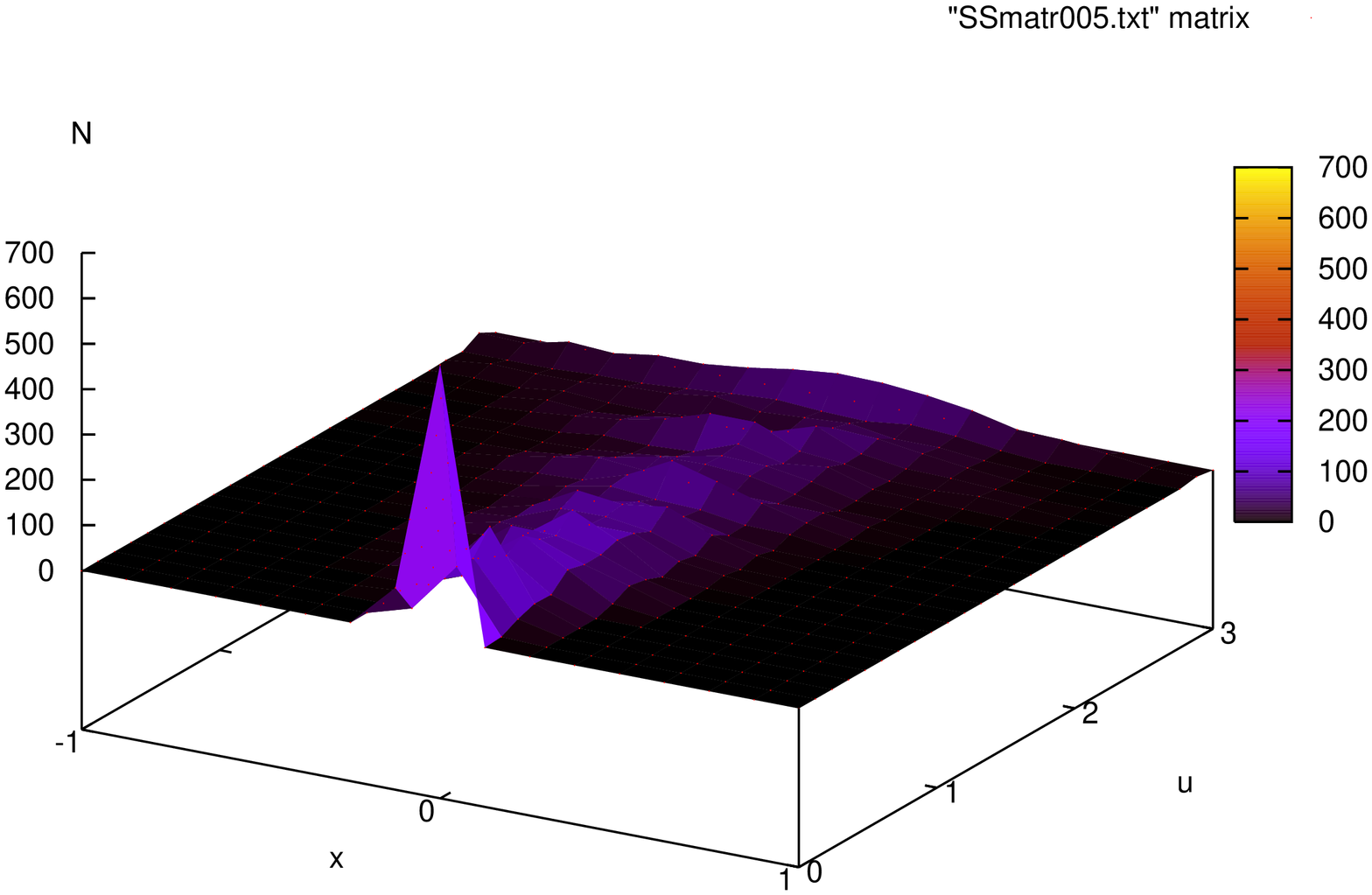, width=.38\textwidth,angle=270}}} \\
  \mbox{\subfigure[$t=20$.]%
{\epsfig{bbllx=50pt,bblly=50pt,bburx=554pt,bbury=770pt,%
figure=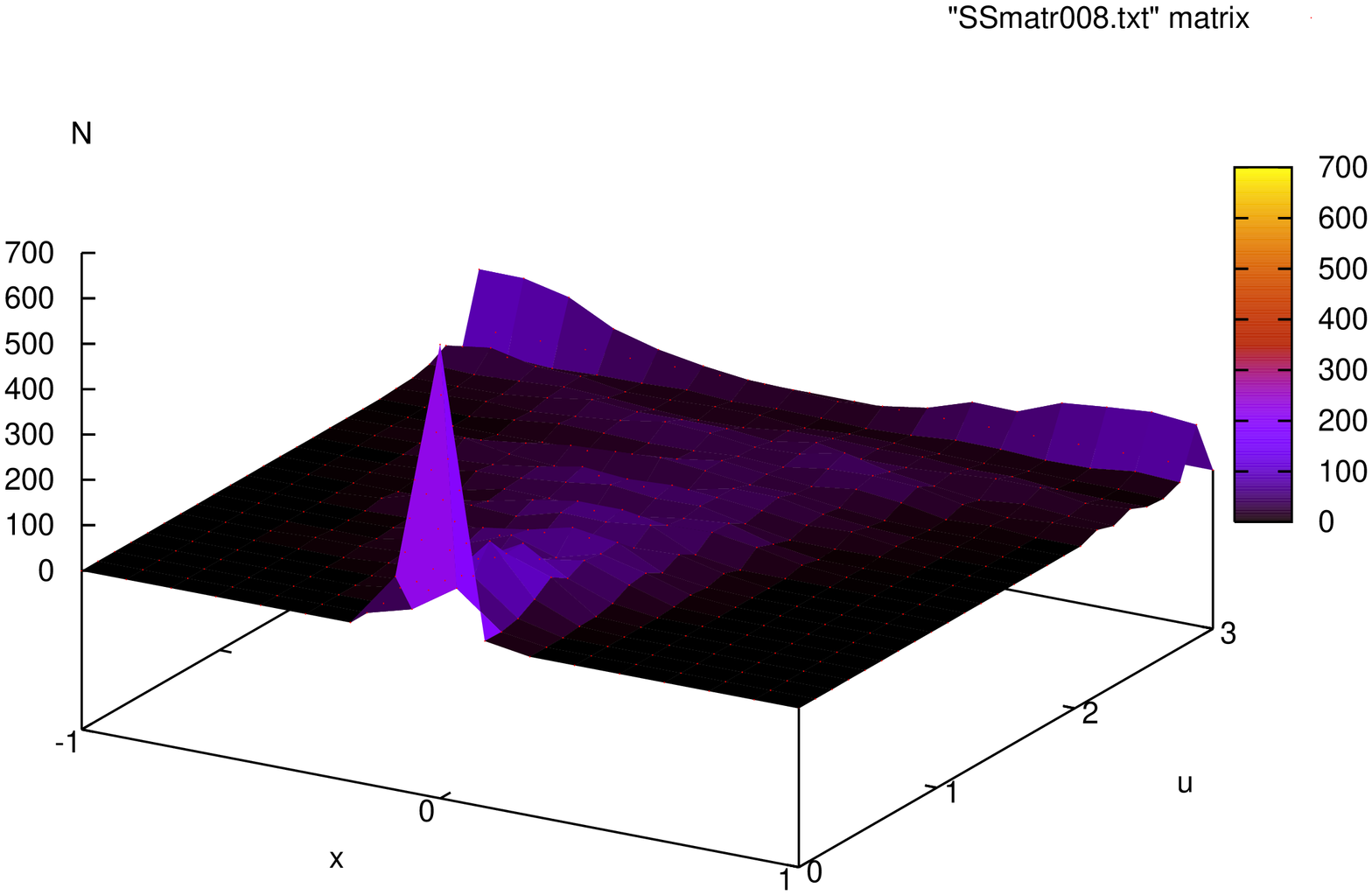, width=.38\textwidth,angle=270}}\quad
    \subfigure[$t=50$.]%
{\epsfig{bbllx=50pt,bblly=50pt,bburx=554pt,bbury=770pt,%
figure=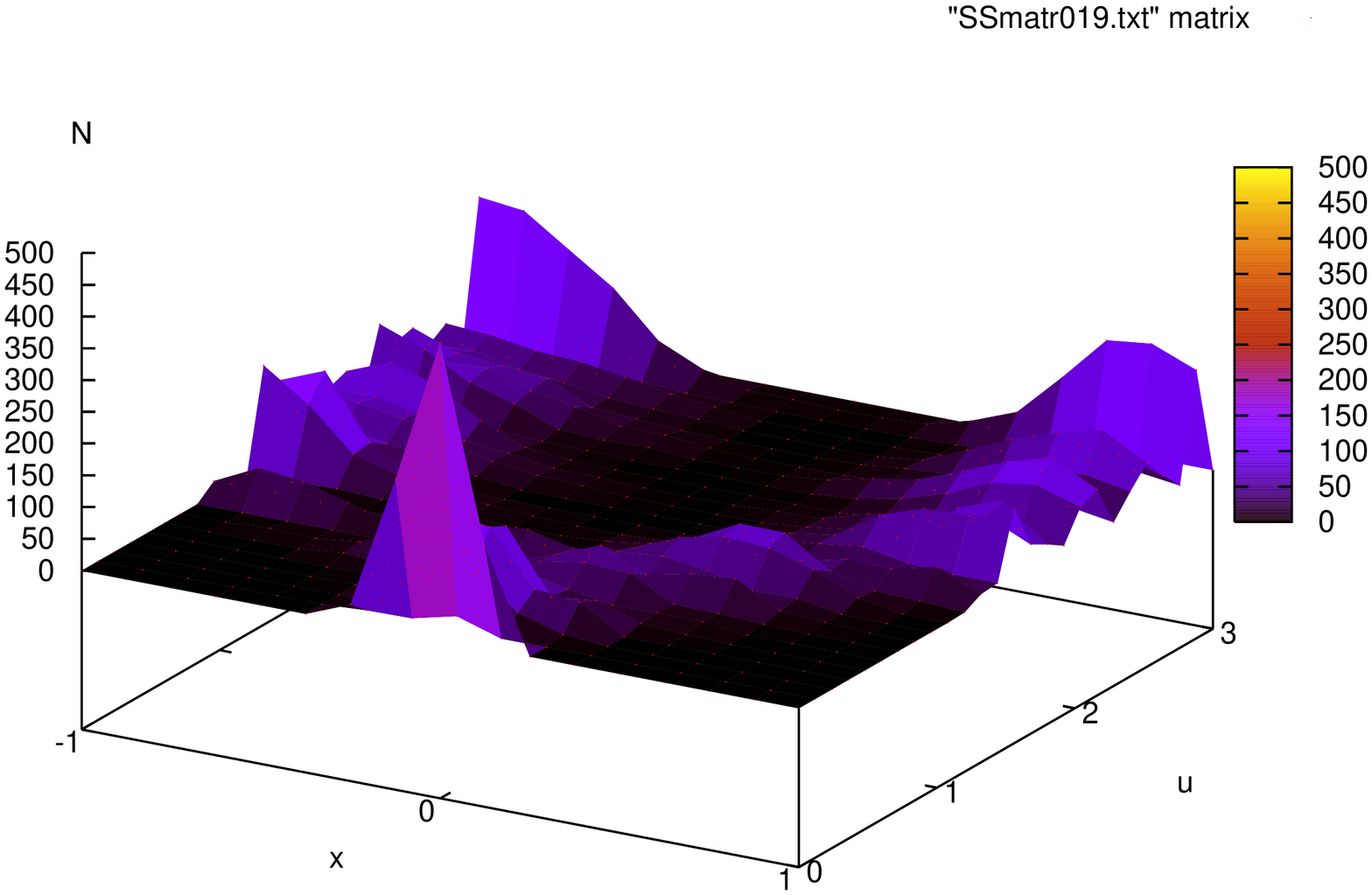, width=.38\textwidth,angle=270}}}
  \caption{Simulations of example 3. The parameters are $N=1000, s=0.03,
  m=0.03$ and $\delta=1$. The initial condition is composed of $N$
  individuals located at 0 and with trait values $3i/N$ for $1\leq
  i\leq N$.}
  \label{fig:ex9}
\end{figure}

When parameters vary, the simulations of this microscopic model
can show different ways of colonization. As an illustration, we
give an example (Fig.~\ref{fig:ex9}) where the interaction range
$\delta$ is bigger. The parameters $N$ and $m$ are chosen such
that two clusters appear for large traits and spread over the
trait space in a short time. The initial condition is the same as
in Fig.~\ref{fig:ex8}. We can observe two branches linking the
initial cluster with the two extreme valued clusters
(Fig.~\ref{fig:ex9}(c,d)).

\bigskip~\bigskip

{\bf Acknowledgments.} The authors would like to thank Laurent
Desvillettes who pointed out to our attention the interest of
combining space and traits, and the article \cite{Pal06}. They
also thank the other participants of the ACI "Structured
Populations" and more specifically R\'egis Ferri\`ere  for
fruitful discussions.

\end{document}